\documentclass[preprint,12pt]{elsarticle}

\newtheorem{theorem}{Theorem}[section]
\newtheorem{corollary}{Corollary}[section]
\newtheorem{lemma}{Lemma}[section]
\newtheorem{remark}{Remark}[section]
\newtheorem{proposition}{Proposition}[section]

\usepackage{amssymb}

\begin{document}

\begin{frontmatter}

\title{Cram\'{e}r large deviation expansions for martingales under Bernstein's condition}
\author{Xiequan Fan$^*$,\ \ \ \ \ \ Ion Grama\ \ \  \ \ and\ \ \ \ \ \ Quansheng Liu}
 \cortext[cor1]{\noindent Corresponding author. \\
\mbox{\ \ \ \ }\textit{E-mail}: fanxiequan@hotmail.com (X. Fan), \ \ \ \ \  ion.grama@univ-ubs.fr (I. Grama),\\
\mbox{\ \ \ \ \ \ \ \ \ \ \  \ \ \ \ }quansheng.liu@univ-ubs.fr (Q. Liu). }
\address{Universit\'{e} de Bretagne-Sud, LMBA, UMR CNRS 6205,
 Campus de Tohannic,\\ 56017 Vannes, France}

\begin{abstract}
An expansion of large deviation probabilities for martingales is given, which extends the classical result due to Cram\'{e}r to the case of martingale differences satisfying the conditional Bernstein condition.
The upper bound of the range of validity and the remainder of our expansion is the same as in the Cram\'{e}r result and therefore are optimal.
Our result implies a moderate deviation principle for martingales.
\end{abstract}

\begin{keyword} expansions of large deviations;
Cram\'{e}r type large deviations; large deviations; moderate deviations; exponential inequality; Bernstein's condition;  central limit theorem
\vspace{0.3cm}
\MSC Primary 60G42; 60F10; 60E15; Secondary 60F05
\end{keyword}

\end{frontmatter}


\section{Introduction}
Consider a sequence of independent and identically distributed (i.i.d.) centered real random variables $\xi_{1},..., \xi_{n}$
satisfying Cram\'{e}r's condition $\mathbb{E}\exp\{ c_{0}|\xi_{1}|\}<\infty,$ for some constant $c_{0}>0.$
Denote $\sigma^2=\mathbb{E}\xi_{1}^2$ and $X_n=\sum_{i=1}^{n}\xi_{i}.$
In 1938, Cram\'{e}r \cite{Cramer38} established an asymptotic expansion of the probabilities of large deviations of $X_n$, based on the powerful technique of conjugate distributions (see also Esscher \cite{Esscher32}).
The results of Cram\'{e}r imply that, uniformly in $1\leq x =  o(n^{1/2}), $
\begin{equation}
\log \frac {\mathbb{P}(X_n> x\sigma\sqrt{n})} {1-\Phi(x)}  =  O \bigg( \frac{x^3}{\sqrt{n}}\bigg) \ \ \mbox{as} \ \ n \rightarrow \infty,
\label{cramer001}
\end{equation}
where $\Phi(x)=\frac{1}{\sqrt{2\pi}}\int_{-\infty}^{x}\exp\{-t^2/2\}dt$ is the standard normal distribution.
Various large deviation expansions for sums of independent random variables have been obtained by many authors, see for instance  Feller \cite{Fl43}, Petrov \cite{Pe54}, Rubin and Sethuraman \cite{RS65}, Statulevi\v{c}ius \cite{S66}, Saulis and Statulevi\v{c}ius \cite{SS78} and Bentkus and Ra\v{c}kauskas  \cite{BR90}.
We refer to the book of Petrov \cite{Petrov75} and the references therein for a detailed account.
Despite the fact that the case of sums of independent random variables is well studied, there are only a few results on expansions of type (\ref{cramer001}) for martingales: see
Bose \cite{Bose86a,Bose86b}, Ra\v{c}kauskas \cite{Rackauskas90,Rackauskas95,Rackauskas97}, Grama \cite{Grama95,G97} and Grama and Haeusler \cite{GH00,GH06}. It is also worth noting that limit theorems for  large and moderate deviation principle for martingales have been proved by several authors, see e.g.\ Liptser and Pukhalskii \cite{LP92}, Gulinsky and Veretennikov \cite{GV93}, Gulinsky, Liptser and Lototskii \cite{GLL94}, Gao \cite{G96}, Dembo \cite{D96},  Worms \cite{W01} and Djellout \cite{D02}. However, these theorems are less precise than large deviation expansions of type (\ref{cramer001}).

Let $(\xi _i,\mathcal{F}_i)_{i=0,...,n}$ be a sequence of  square
integrable martingale differences defined on a probability space $(\Omega, \mathcal{F}, \mathbb{P})$, where $\xi _0=0$ and $\{\emptyset, \Omega\}=\mathcal{F}_0\subseteq ...\subseteq
\mathcal{F}_n\subseteq \mathcal{F}$.
Denote $X_n=\sum_{i=1}^{n}\xi_{i}.$
Assume that there exist absolute constants  $H>0$ and $N\geq 0$ such that $\max_i\left| \xi_i\right| \leq H$ and
$\left| \sum_{i=1}^n \mathbb{E}(\xi _i^2|\mathcal{F}_{i-1})-n\right| \leq N^2.$
Here and hereafter, the equalities and inequalities between random variables are understood in the $\mathbb P$-almost sure sense.
From the results in Grama and  Haeusler \cite{GH00}, it follows that, for any constant $\alpha > 0$ and  $\alpha \sqrt{\log n} \leq x = o\left( n^{1/6}\right),$
\begin{equation}
\frac {\mathbb{P}\left( X_n>x\sqrt{n}\right)} {1-\Phi \left( x\right)} =
1+O\bigg( (H+N)
 \frac {x^3}{\sqrt{n}} \bigg)
\label{cramer002}
\end{equation}
and, for any $0 \leq x = O\left( \sqrt{\log n} \right)$,
\begin{equation}
\frac {\mathbb{P}\left( X_n>x\sqrt{n}\right)} {1-\Phi \left( x\right)} =
1+O\bigg( (H+N) (1+x)\frac {\log n}{\sqrt{n}} \bigg)
\label{cramer003}
\end{equation}
as $n\rightarrow \infty$
(se also \cite{G97,GH06} for more results in the last range).
In this paper we extend the expansions (\ref{cramer002}) and (\ref{cramer003}) to the case of martingale differences $(\xi _i,\mathcal{F}_i)_{i=0,...,n}$ satisfying
the conditional Bernstein condition,
\begin{equation}
|\mathbb{E}(\xi_{i}^{k}  | \mathcal{F}_{i-1})| \leq \frac{1}{2}k!H^{k-2} \mathbb{E}(\xi_{i}^2 | \mathcal{F}_{i-1}),\ \ \ \ \ \mbox{for}\ \ k\geq 3\ \ \mbox{and} \ \ 1\leq i\leq n,
\label{Bernst cond}
\end{equation}
where $H$ is a positive absolute constant.
Note that in the i.i.d.\ case Bernstein's condition (\ref{Bernst cond}) is equivalent to Cram\'{e}r's condition (see Section \ref{seca})
and therefore  (\ref{cramer002}) implies Cram\'{e}r's expansion (\ref{cramer001}).
It is worth stressing that the remainder in expansion (\ref{cramer002}) is of the same order as that in (\ref{cramer001}) in the stated range
and therefore cannot be improved.
As to the remainder in (\ref{cramer003}), from the rate of convergence result in Bolthausen \cite{Bo82} we conclude that it is also optimal.

Another objective of the paper is to find an asymptotic expansion of large deviation for martingales in a wider range than that of (\ref{cramer002}).
From Theorems \ref{th0} and \ref{th1} of the paper it follows that, for any constant $\alpha > 0$ and  $\alpha \sqrt{\log n} \leq x =o\left( n^{1/2}\right),$
\begin{equation}
\log \frac { \mathbb{P}\left( X_n>x\sqrt{n}\right) } { 1-\Phi \left( x\right)} =
 O\bigg(\frac {x^3}{\sqrt{n}} \bigg)  \ \ \mbox{as} \ \ n \rightarrow \infty.
\label{cramer004aaa}
\end{equation}
This improves the corresponding result in \cite{GH00} where (\ref{cramer004aaa}) has been established only in the range $x \in [\alpha \sqrt{\log n}, \alpha_1 n^{1/4}]$ for some absolute constant $\alpha_1>0$.
The upper bound of the range and the remainder in expansion (\ref{cramer004aaa}) cannot be improved since they are of the same order as in the Cram\'{e}r's expansion (\ref{cramer001}).

The idea behind our approach is similar to that of Cram\'{e}r for independent random variables with corresponding adaptations to the martingale case.
We make use of the conjugate multiplicative martingale for changing the probability measure as proposed in Grama and  Haeusler \cite{GH00} (see also \cite{FGL12}).
However, we refine \cite{GH00} in two aspects. First, we relax the boundedness condition $|\xi_i |\le L,$ replacing it by Bernstein's condition (\ref{Bernst cond}).
Secondly, we establish upper and lower bounds for the large deviation probabilities in the range $x \in [0, \alpha_1 n^{1/2})$ thus enlarging the range $x \in [0, \alpha_1 n^{1/4}]$ established in \cite{GH00}.
In the proof we make use of a rate of convergence result for martingales under the conjugate measure. It is established under the Bernstein condition (\ref{Bernst cond}), unlike \cite{GH00} where it is established only for bounded martingale differences. As a consequence, we  improve the result on the rate of convergence in the central limit theorem (CLT) due to Bolthausen \cite{Bo82} (see Theorem \ref{th4} below).

The paper is organized as follows. Our main results are stated and discussedin in Section \ref{sec2}.
A rate of convergence in the CLT for martingales is given in Section \ref{sec2.5}.
Section \ref{sec3} contains auxiliary assertions used in the proofs of the main results.
Proofs are deferred to Sections \ref{sec3.1},  \ref{sec3.2} and \ref{sec4}.
We clarify the relations among the conditions of Bernstein, Cram\'{e}r and Sakhanenko in Section \ref{seca}.

Throughout the paper, $c$ and $c_\alpha,$ probably supplied with some indices,
denote respectively a generic positive absolute constant and a generic positive constant depending only on $\alpha.$
Moreover, $\theta_i$'s stand for values satisfying $\left| \theta_i \right| \leq 1$.

\section{\textbf{Main results} \label{sec2}}
\subsection{Main theorems}
Assume that we are given a sequence of martingale differences $(\xi _i,\mathcal{F}_i)_{i=0,...,n}$, defined on some
 probability space $(\Omega ,\mathcal{F},\mathbb{P})$,  where $\xi
_0=0 $,  $\{\emptyset, \Omega\}=\mathcal{F}_0\subseteq ...\subseteq \mathcal{F}_n\subseteq
\mathcal{F}$ are increasing $\sigma$-fields and $(\xi _i)_{i=1,...,n}$ are allowed to depend on $n$. Set
\begin{equation}
X_{0}=0,\ \ \ \ \ X_k=\sum_{i=1}^k\xi _i,\quad k=1,...,n.  \label{xk}
\end{equation}
Let $\left\langle X\right\rangle $ be the quadratic characteristic of the
martingale $X=(X_k,\mathcal{F}_k)_{k=0,...,n}:$%
\begin{equation}\label{quad}
\left\langle X\right\rangle _0=0,\ \ \ \ \ \left\langle X\right\rangle _k=\sum_{i=1}^k\mathbb{E}(\xi _i^2|\mathcal{F}
_{i-1}),\quad k=1,...,n.
\end{equation}
In the sequel we shall use the following conditions:

\begin{description}
\item[(A1)]  There exists a number $\epsilon \in (0, \frac12]$ such that
\[
|\mathbb{E}(\xi_{i}^{k}  | \mathcal{F}_{i-1})| \leq \frac12 k!\epsilon^{k-2} \mathbb{E}(\xi_{i}^2 | \mathcal{F}_{i-1}),\ \ \ \ \ \textrm{for}\ k\geq 3\ \ \textrm{and}\ \ 1\leq i\leq n;
\]
\item[(A2)]  There exists a number  $ \delta\in [0, \frac12]$ such that
$ \left| \left\langle X\right\rangle _n-1\right| \leq  \delta^2.$
\end{description}

Note that in the case of normalized sums of i.i.d.\ random variables conditions (A1) and (A2) are satisfied with $\epsilon = \frac {1} {\sigma \sqrt n}$ and $\delta = 0$
(see conditions (A1$'$) and (A2$'$) below).
In the case of martingales $\epsilon$ and $\delta$ usually depend on $n$ such that $\epsilon=\epsilon_n \rightarrow 0$ and $\delta=\delta_n \rightarrow 0$.

The following two theorems give upper and lower bounds for large deviation probabilities.

\begin{theorem}\label{th0}
Assume conditions (A1) and (A2). Then, for any constant $\alpha \in (0,1)$ and all $0\leq x \leq \alpha \, \epsilon^{-1},$ we have
\begin{equation}\label{t0ie1}
 \frac{\mathbb{P}(X_n>x)}{1-\Phi \left( x\right)}\leq \exp \Bigg\{ c_{\alpha} \bigg( x^3 \epsilon + x^2 \delta^2 \bigg)  \Bigg\}
\Bigg(\frac{}{} 1+c_{\alpha}\,  (1+ x)\left( \epsilon \left| \log  \epsilon
 \right| + \delta \right) \Bigg)
\end{equation}
and
\begin{equation}
 \frac{\mathbb{P}(X_n<-x)}{\Phi \left( -x\right)}\leq \exp\Bigg\{ c_{\alpha} \bigg( x^3 \epsilon + x^2 \delta^2 \bigg)  \Bigg\}
\Bigg(\frac{}{} 1+c_{\alpha}\, (1+ x)\left( \epsilon \left| \log  \epsilon
 \right| + \delta \right) \Bigg),
\end{equation}
where the constant $c_{\alpha}$ does not depend on $(\xi _i,\mathcal{F}_i)_{i=0,...,n}$, $n$ and $x$.
\end{theorem}

\begin{theorem}\label{th1}
Assume conditions (A1) and (A2). Then there is an absolute constant $\alpha_0 >0$ such that,
for all $0\leq x \leq \alpha_0 \,\epsilon^{-1}$ and $\delta \leq \alpha_0$,
\begin{equation} \label{t1ie1}
\frac{\mathbb{P}(X_n>x)}{1-\Phi \left( x\right)}\geq \exp\Bigg\{-c_{\alpha_0} \! \bigg( x^3 \epsilon  + x^2 \delta^2 + (1+ x) \left( \epsilon \left| \log  \epsilon
 \right|+  \delta \right)\bigg) \Bigg\}
\end{equation}
and
\begin{equation}
  \frac{\mathbb{P}(X_n<-x)}{\Phi \left( -x\right)} \geq  \exp\Bigg\{-c_{\alpha_0} \! \bigg( x^3 \epsilon  + x^2 \delta^2 + (1+ x) \left( \epsilon \left| \log  \epsilon\right|+  \delta \right)\bigg) \Bigg\},
\end{equation}
where the constants $\alpha_0$ and $c_{\alpha_0}$ do not depend on $(\xi _i,\mathcal{F}_i)_{i=0,...,n}$, $n$ and $x$.
\end{theorem}

Using the inequality $|e^x -1 | \leq   e^{ \alpha } |x|$ valid for $|x| \leq \alpha, $ from Theorems \ref{th0} and \ref{th1},
we obtain the following improvement of the main result of \cite{GH00}.
\begin{corollary}\label{corollary01}
Assume conditions (A1) and (A2). Then there is an absolute constant $\alpha_0 >0$ such that,
for all $0\leq x \leq \alpha_0 \, \min\{ (\epsilon \left| \log  \epsilon \right|)^{-1}, \delta^{-1}\},$
\begin{equation}\label{fth1}
\frac{\mathbb{P}(X_n>x)}{1-\Phi \left( x\right)}=\exp\{\theta_{1}c_{\alpha_0} x^3 \epsilon \}
\left(\frac{}{} 1+\theta_{2}c_{\alpha_0}(1+x)( \epsilon \left| \log  \epsilon
 \right| + \delta ) \right)
\end{equation}
and
\begin{equation}\label{fth2}
\frac{\mathbb{P}(X_n<-x)}{\Phi \left( -x\right) } = \exp  \{\theta_{3}c_{\alpha_0}  x^3 \epsilon \}
\left(\frac{}{} 1+\theta_{4}c_{\alpha_0}(1+x)( \epsilon \left| \log \epsilon
 \right| + \delta ) \right),
\end{equation}
where $c_{\alpha_0}$ does not depend on $n, x$ but $ \theta_i $ possibly depend on $(\xi _i,\mathcal{F}_i)_{i=0,...,n}$, $n$ and $x.$
\end{corollary}

For bounded martingale differences $|\xi_i|\leq \epsilon$ under condition (A2),  Grama and  Haeusler \cite{GH00} proved the asymptotic expansions  (\ref{fth1}) and (\ref{fth2}) for $x \in [0,  \alpha_1  \min\{\epsilon^{-1/2}, \delta^{-1}\}]$ and some small absolute constant $\alpha_1 \in (0, \frac18]$. Thus Corollary \ref{corollary01} extends the asymptotic expansions of \cite{GH00} to a  larger range $x \in [0,  \alpha_0 \, \min\{ (\epsilon \left| \log  \epsilon\right|)^{-1}, \delta^{-1}\})$ and non bounded martingale differences.

\subsection{Remarks on the main theorems}
Combining the inequalities (\ref{t0ie1}) and (\ref{t1ie1}), we conclude that under (A1) and (A2) there is an absolute constant $\alpha_0 >0$ such that,
for all $0\leq x \leq \alpha_0 \,\epsilon^{-1}$ and $\delta \leq \alpha_0$,
\begin{equation}\label{fsvdf1}
\Bigg| \log \frac{\mathbb{P}(X_n>x)}{1-\Phi \left( x\right)} \Bigg|  \leq   c_{\alpha_0} \bigg( x^3 \epsilon + x^2 \delta^2+ (1+x) \left(\epsilon \left| \log  \epsilon\right|+  \delta\right) \bigg).
\end{equation}
We show that this result can be regarded as a refinement of the moderate deviation principle (MDP) in the framework where (A1) and (A2) hold.
Assume that (A1) and (A2) are satisfied with $\epsilon=\epsilon_n\rightarrow 0$ and $\delta=\delta_n\rightarrow 0$ as $n\rightarrow \infty$.
Let $a_n$ be any sequence of real numbers satisfying $a_n \rightarrow \infty$ and $a_n\epsilon_n\rightarrow 0$
as $n\rightarrow \infty$.
Then inequality (\ref{fsvdf1}) implies the MDP for $(X_n)_{n\geq 1}$ with the speed $a_n$ and rate function $x^2/2.$
Indeed, using the inequalities
\[
\frac{1}{\sqrt{2 \pi}(1+x)} e^{-x^2/2} \leq 1-\Phi ( x ) \leq \frac{1}{\sqrt{ \pi}(1+x)} e^{-x^2/2}, \ \  x\geq 0,
\]
we deduce that, for any $x\geq0$,
\begin{eqnarray*}
\lim_{n\rightarrow \infty}\frac{1}{a_n^2}\log \mathbb{P}(X_n>a_nx) = - \frac{x^2}{2} .
\end{eqnarray*}
By a similar argument, we also have, for any $x\geq0$,
\begin{eqnarray*}
\lim_{n\rightarrow \infty}\frac{1}{a_n^2}\log \mathbb{P}(X_n< -a_nx) = - \frac{x^2}{2}.
\end{eqnarray*}
The last two equalities are equivalent to the statement that: for each Borel set $B$,
\begin{eqnarray*}
- \inf_{x \in B^o}\frac{x^2}{2} &\leq & \liminf_{n\rightarrow \infty}\frac{1}{a_n^2}\log \mathbb{P}\left(\frac{1}{a_n} X_n \in B \right) \\
 &\leq& \limsup_{n\rightarrow \infty}\frac{1}{a_n^2}\log \mathbb{P}\left(\frac{1}{a_n} X_n \in B \right) \leq  - \inf_{x \in \overline{B}}\frac{x^2}{2} \, ,
\end{eqnarray*}
where $B^o$ and $\overline{B}$ denote the interior and the closure of $B$ respectively, see Lemma 4.4 of \cite{HT}.
Similar results can be found in Gao \cite{G96} for the martingale differences satisfying the conditional Cram\'{e}r condition  $||\mathbb{E}(\exp\{ c_{0}|\xi_{i}|\}| \mathcal{F}_{i-1})||_\infty<\infty$.

To show that our results are sharp, assume that $\xi_{i}=\eta_{i}/\sqrt{n}$, where $(\eta _i,\mathcal{F}_i)_{i=1,...,n}$ is a sequence of martingale differences satisfying the following conditions:
\begin{description}
\item[(A1$'$)] (\emph{Bernstein's condition}) There exists a positive absolute constant $H$ such that
\[
|\mathbb{E}(\eta_{i}^{k}  | \mathcal{F}_{i-1})| \leq \frac12 k!H^{k-2} \mathbb{E}(\eta_{i}^2 | \mathcal{F}_{i-1}), \    \ \ \textrm{for}\  k\geq 3\ \ \textrm{and}\ \ 1\leq i\leq n;
\]
\item[(A2$'$)] There exists an absolute constant $N\geq 0$ such that $\left|\sum_{i=1}^{n}\mathbb{E}(\eta^2_{i} | \mathcal{F}_{i-1}) -n  \right| \leq  N^2$.
\end{description}
These conditions are satisfied with some $H>0$ and $N=0$ if, for instance, $\eta_{1}, \eta_{2},..., \eta_{n}$ are i.i.d.\ random variables with  finite exponential moments (see Section \ref{seca} for an explicit expression of the positive absolute constant $H$).

\begin{corollary}\label{co13}
Assume  conditions (A1$'$) and (A2\,$'$). Then there is an absolute constant $\alpha_{2}>0$ such that for any absolute constant $\alpha_{1}>0$ and
all $\alpha_{1}\sqrt{\log n} \leq x \leq \alpha_{2} n^{1/2},$ we have
\begin{equation}\label{co13ie1}
\log \frac{\mathbb{P}( \sum_{i=1}^{n}\eta_{i}>x\sqrt{n})}{1-\Phi \left( x\right)}=    O\left( (H+N) \frac{x^3}{\sqrt{n}} \right)
\end{equation}
and
\begin{equation}\label{co13ie2}
\log \frac{\mathbb{P}( \sum_{i=1}^{n}\eta_{i}<- x\sqrt{n})}{\Phi \left( -x\right) }=   O\left( (H+N) \frac{x^3}{\sqrt{n}} \right)
\end{equation}
as $n\rightarrow \infty$.
\end{corollary}

It is worth noting that the remainders of the expansions (\ref{co13ie1}) and (\ref{co13ie2}) are of the
same order as in (\ref{cramer001}) and therefore are optimal.

\begin{corollary}\label{co1}
Assume conditions (A1$'$) and (A2$'$). Then, for all $0\leq x = O\left( \sqrt{\log n} \right)$,
\begin{equation}\label{corie1}
\frac{\mathbb{P}(\sum_{i=1}^{n}\eta_{i}>x\sqrt{n})}{1-\Phi \left( x\right)}= 1 + O\left( (H+N) (1+x) \frac{\log n}{\sqrt{n}}   \right)
\end{equation}
and
\begin{equation}\label{corie2}
\frac{\mathbb{P}(\sum_{i=1}^{n}\eta_{i}<-x\sqrt{n})}{\Phi \left( -x\right) }= 1 + O \left( (H+N) (1+x) \frac{\log n}{\sqrt{n}}   \right)
\end{equation}
as $n\rightarrow \infty$.
\end{corollary}

Notice that (\ref{corie1}) extends expansion (\ref{cramer003}) proved in Grama and  Haeusler \cite{GH00} to the case of martingale differences satisfying the conditional Bernstein condition (A1$'$). The Remark 2.1 of \cite{GH00} and the sharp rate of convergence in the CLT due to Bolthausen \cite{Bo82} hint that the remainders of the expansions (\ref{corie1}) and (\ref{corie2}) are sharp.

\begin{corollary}\label{co12}
Assume  conditions (A1$'$) and (A2\,$'$). Then, for any  absolute constant $\alpha > 0$ and  $\alpha \sqrt{\log n} \leq x =o\left( n^{1/6}\right),$
\begin{equation}\label{co12ie1}
\frac{\mathbb{P}(\sum_{i=1}^{n}\eta_{i}>x\sqrt{n})}{1-\Phi \left( x\right)}= 1 +  O\bigg( (H+N)\frac{x^3}{\sqrt{n}}\bigg)
\end{equation}
and
\begin{equation}\label{co12ie2}
\frac{\mathbb{P}( \sum_{i=1}^{n}\eta_{i}<- x\sqrt{n})}{\Phi \left( -x\right) }=1 +  O\bigg( (H+N)\frac{x^3}{\sqrt{n}}\bigg)
\end{equation}
as $n\rightarrow \infty$.
\end{corollary}

The remainders of the expansions (\ref{co12ie1}) and (\ref{co12ie2}) are of the
same order as in (\ref{cramer001}) in the stated range and therefore cannot be improved.

\begin{remark}
The results formulated above are proved under Bernstein's condition (A1$'$). But they are also valid under some equivalent conditions which are stated in Section \ref{seca}.
\end{remark}

\section{\textbf{Rates of convergence in the CLT}\label{sec2.5}}
Let $(\xi _i,\mathcal{F}_i)_{i=0,...,n}$ be a sequence of martingale differences
satisfying condition (A1) and $X=(X_k,\mathcal{F}_k)_{k=0,...,n}$
be the corresponding martingale defined by (\ref{xk}). For any real $\lambda$ satisfying $|\lambda| < \epsilon^{-1} ,$
consider the \emph{exponential multiplicative martingale} $Z(\lambda
)=(Z_k(\lambda ),\mathcal{F}_k)_{k=0,...,n},$ where
\[
Z_k(\lambda )=\prod_{i=1}^k\frac{e^{\lambda \xi _i}}{\mathbb{E}(e^{\lambda \xi _i}|
\mathcal{F}_{i-1})},\quad k=1,...,n,\quad Z_0(\lambda )=1.  \label{C-1}
\]
For each $k=1,...,n,$ the random variable $Z_k(\lambda
) $ defines a probability density on $(\Omega ,\mathcal{F},\mathbb{P}).$ This allows us to introduce, for $|\lambda|
 <\epsilon^{-1},$ the \emph{conjugate probability measure} $\mathbb{P}_\lambda $ on $(\Omega ,%
\mathcal{F})$ defined by
\begin{equation}
d\mathbb{P}_\lambda =Z_n(\lambda )d\mathbb{P}.  \label{f21}
\end{equation}
Denote by $\mathbb{E}_{\lambda}$ the expectation with respect to $\mathbb{P}_{\lambda}$.
For all $i=1,\dots,n$, let
\[
\eta_i(\lambda)=\xi_i - b_i(\lambda)\ \ \ \ \ \ \ \ \  \textrm{and} \ \ \ \  \ \ \  \ \  b_i(\lambda)=\mathbb{E}_{\lambda}(\xi_i |\mathcal{F}_{i-1}).  \]
We thus obtain the well-known semimartingale decomposition:
\begin{equation}
X_k=Y_k(\lambda )+B_k(\lambda ),\quad\quad\quad k=1,...,n, \label{xb}
\end{equation}
where $Y(\lambda )=(Y_k(\lambda ),\mathcal{F}_k)_{k=1,...,n}$ is the \emph{%
conjugate martingale} defined as
\begin{equation}\label{f23}
Y_k(\lambda )=\sum_{i=1}^k\eta _i(\lambda ),\quad\quad\quad k=1,...,n,
\end{equation}
and
$B(\lambda )=(B_k(\lambda ),\mathcal{F}_k)_{k=1,...,n}$ is the \emph{%
drift process} defined as
\[
B_k(\lambda )=\sum_{i=1}^kb_i(\lambda ),\quad\quad\quad k=1,...,n.
\]
In the proofs of Theorems \ref{th0} and \ref{th1}, we make use of the following assertion, which gives
us a rate of convergence in the central limit theorem for the conjugate
martingale $Y( \lambda  )$ under the probability measure $\mathbb{P}_{ \lambda  }.$
\begin{lemma}
\label{LEMMA4}
Assume conditions (A1) and (A2). Then, for all $0\leq \lambda < \epsilon ^{-1}$,
\[
\sup_{x}\left| \mathbb{P}_\lambda (\, Y_n(\lambda )\leq x)-\Phi (x)\right| \leq
c\left( \lambda \, \epsilon +\epsilon \left| \log   \epsilon \right| +\delta \right) .
\]

\end{lemma}

If $\lambda=0$, then
$Y_n(\lambda )=X_{n}$ and $\mathbb{P}_\lambda=\mathbb{P}.$
So Lemma \ref{LEMMA4} implies the following theorem.
\begin{theorem}\label{th4}
Assume conditions (A1) and (A2). Then
\begin{equation}\label{f30}
\sup_{x}|\mathbb{P}(X_n \leq x)-\Phi \left( x\right)| \leq c \,(\epsilon \left| \log  \epsilon\right|+\delta ).
\end{equation}
\end{theorem}
\begin{remark}
By inspecting the proof of Lemma \ref{LEMMA4}, we can see that Theorem \ref{th4} holds true when
 condition (A1) is replaced by the following weaker one:
\begin{description}
\item[(C1)]  There exists a number  $ \epsilon \in (0, \frac 1 2]$ depending on $n$ such that
\[
|\mathbb{E}(\xi_{i}^{k}  | \mathcal{F}_{i-1})| \leq  \epsilon^{k-2} \mathbb{E}(\xi_{i}^2 | \mathcal{F}_{i-1}), \ \ \ \textrm{for}\ k=3,5 \ \textrm{and} \ 1\leq i\leq n.
\]
\end{description}
\end{remark}
\begin{remark}
Bolthausen (see Theorem 2 of \cite{Bo82}) showed that if $|\xi_{i}|\leq \epsilon$ and condition (A2) holds, then
\begin{equation}\label{f31}
\sup_{x }|\mathbb{P}(X_n \leq x)-\Phi \left( x\right)| \leq c_1 \,(\epsilon^3 n  \log n  +\delta).
\end{equation}
We note that Theorem  \ref{th4} implies Bolthausen's inequality (\ref{f31}) under the less restrictive condition (A1).
Indeed, by condition (A2), we have $3/4 \leq \langle X \rangle_n \leq n \epsilon ^2$ and then $\epsilon \geq \sqrt{3/(4n)}$. For $\epsilon \leq 1/2$, it is easy to see that $\epsilon^3 n   \log n  \geq 3\,\epsilon|\log\epsilon | /4$.
Thus, inequality (\ref{f30}) implies (\ref{f31}) with $c_1=4c/3$.
\end{remark}

\section{\textbf{Auxiliary results}\label{sec3}}
In this section, we establish some auxiliary lemmas which will be used in the proofs of Theorems \ref{th0} and \ref{th1}.
We first prove upper bounds for the conditional moments.
\begin{lemma}
\label{l11} Assume condition (A1). Then
\[
|\mathbb{E}(\xi_{i}^k | \mathcal{F}_{i-1})| \leq   6 k! \epsilon^{k}, \ \ \ \textrm{for} \  k\geq 2,
\]
and
\[
\mathbb{E}(|\xi_{i}|^k | \mathcal{F}_{i-1}) \leq   k! \epsilon^{k-2} \mathbb{E}(\xi_{i}^2 | \mathcal{F}_{i-1}), \ \ \ \textrm{for} \  k\geq 2.
\]
\end{lemma}

\noindent\textbf{Proof.}
By Jensen's inequality and  condition (A1),
\begin{eqnarray*}
\mathbb{E}(\xi_{i}^2 | \mathcal{F}_{i-1})^{2}  \leq   \mathbb{E}(\xi_{i}^4 | \mathcal{F}_{i-1})
 \leq  12 \epsilon^2 \mathbb{E}(\xi_{i}^2|\mathcal{F}_{i-1}),
\end{eqnarray*}
from which we get
\[
\mathbb{E}(\xi_{i}^2 | \mathcal{F}_{i-1}) \leq 12 \epsilon^2 .
\]
We obtain the first assertion. Again by  condition (A1), for $k\geq 3,$
\begin{eqnarray*}
|\mathbb{E}(\xi_{i}^k | \mathcal{F}_{i-1})|
 \leq   \frac12 k!\epsilon^{k-2} \mathbb{E}(\xi_{i}^2 | \mathcal{F}_{i-1})
 \leq  6 k!\epsilon^{k}.
\end{eqnarray*}
If $k$ is even, the second assertion holds obviously. If $k=2l+1$, $l\geq1$, is odd,
by H\"{o}lder's inequality and  condition (A1), it follows that
\begin{eqnarray*}
\mathbb{E}\left(|\xi_{i}|^{2l+1} | \mathcal{F}_{i-1} \right)
&\leq & \mathbb{E}\left(|\xi_{i}|^{l}|\xi_{i}|^{l+1} | \mathcal{F}_{i-1} \right) \leq \sqrt{\mathbb{E}\left( \xi_{i}^{2l} | \mathcal{F}_{i-1} \right)\mathbb{E}\left(\xi_{i}^{2(l+1)} | \mathcal{F}_{i-1} \right) }\\
&\leq&\frac{1}{2}\sqrt{(2l)!(2l+2)!} \epsilon^{2l-1} \mathbb{E}(\xi_{i}^2 | \mathcal{F}_{i-1}) \\
&\leq & (2l+1)!  \epsilon^{2l-1} \mathbb{E}(\xi_{i}^2 | \mathcal{F}_{i-1}) .
\end{eqnarray*}
This completes the proof of Lemma \ref{l11}. \hfill\qed

The following lemma establishes a two sided bound for the drift process $B_n(\lambda ).$
\begin{lemma}
\label{LEMMA-1-1} Assume conditions (A1) and (A2). Then for any constant $\alpha \in (0,1)$ and all $0 \leq \lambda \leq \alpha\, \epsilon^{-1} ,$
\begin{eqnarray}\label{f25sa}
 | B_n(\lambda ) - \lambda |\leq
 \lambda \delta^{2} + c_\alpha \lambda^{2}\epsilon.
\end{eqnarray}
\end{lemma}

\noindent\textbf{Proof.} By the relation between $\mathbb{E}$ and $\mathbb{E}_{\lambda}$ on $\mathcal{F}_i,$
we have
\[
b_i(\lambda )=\frac{\mathbb{E}(\xi
_ie^{\lambda \xi _i}|\mathcal{F}_{i-1})}{\mathbb{E}(e^{\lambda \xi _i}|\mathcal{F}%
_{i-1})},\quad\ \ \ \ \ \   i=1,...,n.
\]
Jensen's inequality and $\mathbb{E}(\xi _i|\mathcal{F}_{i-1})=0$ imply that $\mathbb{E}(e^{\lambda \xi _i}|\mathcal{F}
_{i-1})\geq 1.$ Since
\[
\mathbb{E}(\xi_{i} e^{\lambda\xi_{i}} |\mathcal{F}_{i-1})=\mathbb{E}\left(\xi_{i}(e^{\lambda\xi_{i}}-1)|\mathcal{F}_{i-1} \right)\geq 0,\ \ \ \ \mbox{for}\ \lambda\geq0,
\]
 by Taylor's expansion for $e^x$, we find that
\begin{eqnarray}
B_n(\lambda ) & \leq & \sum_{i=1}^{n}\mathbb{E}(\xi_{i} e^{\lambda \xi_{i}} | \mathcal{F}_{i-1})\nonumber\\
& = & \sum_{i=1}^{n}\mathbb{E}\left(\xi_{i}(e^{\lambda \xi_{i}}-1)| \mathcal{F}_{i-1} \right)\nonumber
\\ & = &  \lambda\langle X\rangle_{n}+ \sum_{i=1}^{n}\sum_{k=2}^{+\infty}\mathbb{E} \left(\frac{\xi_{i}(\lambda\xi_{i})^{k}}{k !}  \Bigg |  \mathcal{F}_{i-1}   \right)   .\label{f26}
\end{eqnarray}
Using condition (A1), we obtain, for any  constant $\alpha \in (0,1)$ and all $0 \leq \lambda \leq \alpha\, \epsilon^{-1} $,
\begin{eqnarray}
  \sum_{i=1}^{n}\sum_{k=2}^{+\infty} \left(\frac{\xi_{i}(\lambda\xi_{i})^{k}}{k !}  \Bigg |  \mathcal{F}_{i-1}   \right)
  & \leq & \sum_{i=1}^{n}\sum_{k=2}^{+\infty}|\mathbb{E}\left( \xi_{i}^{k+1}| \mathcal{F}_{i-1} \right)| \frac{\lambda^k}{k !} \nonumber \\
  & \leq & \frac12\, \lambda^2 \epsilon  \langle X \rangle_{n} \sum_{k=2}^{+\infty}(k+1)(\lambda\epsilon)^{k-2} \nonumber \\
& \leq &c_\alpha \, \lambda^2 \epsilon  \langle X \rangle_{n} .\label{f52}
\end{eqnarray}
Using condition (A2), we get $\langle X \rangle_{n}  \leq 2$ and, for any  constant $\alpha \in (0,1)$ and all $0 \leq \lambda \leq \alpha\, \epsilon^{-1} $,
\begin{eqnarray} \label{f553}
\sum_{i=1}^{n}\sum_{k=2}^{+\infty}\left|\mathbb{E}\left(\frac{\xi_{i}(\lambda \xi_{i})^{k}}{k !} \Bigg| \mathcal{F}_{i-1}\right) \right| \leq 2\, c_\alpha \lambda^2 \epsilon.
\end{eqnarray}
 Condition (A2) together with (\ref{f26}) and (\ref{f553}) imply the upper bound of $B_n(\lambda )$: for any  constant $\alpha \in (0,1)$ and all $0 \leq \lambda \leq \alpha\, \epsilon^{-1} $,
\[
B_n(\lambda ) \leq
\lambda +  \lambda \delta^{2} + 2\,c_\alpha\, \lambda^{2}\epsilon   .
\]
Using  Lemma \ref{l11}, we have, for any  constant $\alpha \in (0,1)$ and all $0 \leq \lambda \leq \alpha\, \epsilon^{-1} $,
\begin{eqnarray}
\mathbb{E}\left(e^{\lambda \xi_{i}} | \mathcal{F}_{i-1} \right) & \leq & 1 + \sum_{k=2}^{+\infty}\left|\mathbb{E} \left(\frac{(\lambda\xi_{i})^k}{k !} \Bigg| \mathcal{F}_{i-1} \right)\right|  \nonumber \\
& \leq & 1+  6\sum_{k=2}^{+\infty} \left(\lambda \epsilon\right)^{k}  \nonumber \\
& \leq & 1+ c_{1,\alpha}\, (\lambda\epsilon )^2 .  \label{sbopdf}
\end{eqnarray}
This inequality together with condition (A2) and (\ref{f553}) imply the lower bound of $B_{n}(\lambda)$: for any  constant $\alpha \in (0,1)$ and all $0 \leq \lambda \leq \alpha\, \epsilon^{-1} $,
\begin{eqnarray*}
  B_n(\lambda ) & \geq & \left(\sum_{i=1}^{n}\mathbb{E}(\xi_{i} e^{\lambda \xi_{i}} | \mathcal{F}_{i-1})\right)\Bigg( 1 + c_{1,\alpha}\, (\lambda\epsilon)^2 \Bigg)^{-1}\nonumber\\
  & \geq &  \left(\lambda\langle X\rangle_{n} -  \sum_{i=1}^{n}\sum_{k=2}^{+\infty}\left|\mathbb{E}\left(\frac{\xi_{i}(\lambda\xi_{i})^{k}}{k !} \Bigg | \mathcal{F}_{i-1} \right) \right| \right)\Bigg( 1 + c_{1,\alpha}\, (\lambda\epsilon)^2 \Bigg)^{-1} \\
  & \geq &  \bigg(\lambda -\lambda \delta^2 -  2\, c_{\alpha} \lambda^2 \epsilon \bigg)\bigg( 1 + c_{1,\alpha}\,(\lambda\epsilon)^2 \bigg)^{-1} \\
  & \geq &  \lambda - \lambda \delta^{2} - (2\, c_{\alpha}+\alpha \, c_{1,\alpha})\, \lambda^{2}\epsilon,
\end{eqnarray*}
where the last line follows from the following inequality, for any  constant $\alpha \in (0,1)$ and all $0 \leq \lambda \leq \alpha\, \epsilon^{-1} $,
\begin{eqnarray*}
\lambda -\lambda \delta^2 -  2\,c_{\alpha} \lambda^2 \epsilon &\geq& \lambda -\lambda \delta^2   - (2\,c_{\alpha}+\alpha \, c_{1,\alpha}) \lambda^2 \epsilon +c_{1,\alpha}\lambda^3 \epsilon^2 \\
&\geq& \bigg(\lambda - \lambda \delta^{2} - (2\, c_{\alpha}+\alpha \, c_{1,\alpha})\, \lambda^{2}\epsilon \bigg)\bigg(1 + c_{1,\alpha} (\lambda\epsilon)^2 \bigg).
\end{eqnarray*}
The proof of Lemma \ref{LEMMA-1-1} is finished.\hfill\qed

Now, consider the predictable cumulant process $\Psi (\lambda )=(\Psi
_k(\lambda ),\mathcal{F}_k)_{k=0,...,n}$ related with the martingale $X$ as follows:
\begin{equation}
\Psi _k ( \lambda )=\sum_{i=1}^k\log \mathbb{E}\left( e^{\lambda \xi _i}|\mathcal{F}_{i-1} \right). \label{f29}
\end{equation}
We establish a two sided bound for the process $\Psi(\lambda ).$
\begin{lemma}
\label{LEMMA-1-2} Assume conditions (A1) and (A2). Then, for any constant $\alpha \in (0,1)$ and all $0 \leq \lambda \leq \alpha\, \epsilon^{-1} ,$
\[
\left| \Psi _n(\lambda )-\frac{\lambda ^2}2\right| \leq   c_\alpha \lambda ^3\epsilon +\frac{\lambda^2 \delta^2}{2}.
\]
\end{lemma}

\noindent\textbf{Proof.}
Since $\mathbb{E}(\xi _i|\mathcal{F}_{i-1})=0$, it is easy to see that
\[
\Psi _n(\lambda )=\sum_{i=1}^n\left( \log \mathbb{E}(e^{\lambda \xi _i}|\mathcal{F}%
_{i-1})-\lambda \mathbb{E}(\xi _i|\mathcal{F}_{i-1})-\frac{\lambda ^2}2\mathbb{E}(\xi _i^2|%
\mathcal{F}_{i-1})\right) +\frac{\lambda ^2}2\left\langle X\right\rangle _n.
\]
Using a two-term  Taylor's expansion of $\log(1+x), x\geq0$,  we obtain
\begin{eqnarray*}
 \Psi _n(\lambda ) - \frac{\lambda ^2}2\left\langle X\right\rangle _n
&=& \sum_{i=1}^n\left( \mathbb{E}(e^{\lambda \xi _i}|\mathcal{F}%
_{i-1})-1-\lambda \mathbb{E}(\xi _i|\mathcal{F}_{i-1})-\frac{\lambda ^2}2\mathbb{E}(\xi _i^2|%
\mathcal{F}_{i-1}) \right) \\
&& - \frac1{2 \, \bigg(1+|\theta|\left(\mathbb{E}(e^{\lambda \xi _i}|\mathcal{F}%
_{i-1})-1 \right) \bigg)^2} \sum_{i=1}^n \left(\frac{}{} \mathbb{E}(e^{\lambda \xi _i}|\mathcal{F}%
_{i-1})-1 \right)^2.
\end{eqnarray*}
Since $\mathbb{E}(e^{\lambda \xi _i}|\mathcal{F}_{i-1})\geq1$, we find that
\begin{eqnarray*}
 \left|\Psi _n(\lambda ) - \frac{\lambda ^2}2\left\langle X\right\rangle _n \right|
&\leq& \sum_{i=1}^n\left| \mathbb{E}(e^{\lambda \xi _i}|\mathcal{F}%
_{i-1})-1-\lambda \mathbb{E}(\xi _i|\mathcal{F}_{i-1})-\frac{\lambda ^2}2\mathbb{E}(\xi _i^2|%
\mathcal{F}_{i-1}) \right|\\
&& + \frac12 \sum_{i=1}^n \left(\frac{}{} \mathbb{E}(e^{\lambda \xi _i}|\mathcal{F}%
_{i-1})-1 \right)^2\\
&\leq& \sum_{i=1}^{n}\sum_{k=3}^{+\infty}\frac{\lambda^{k}}{k !}|\mathbb{E}(\xi_{i}^{k} |\mathcal{F}_{i-1})| + \frac12\sum_{i=1}^{n} \left(\sum_{k=2}^{+\infty}\frac{\lambda^k}{k !}|\mathbb{E}(\xi_{i}^k|\mathcal{F}_{i-1})|\right)^2.
\end{eqnarray*}
In the same way as in the proof of (\ref{f52}), using condition (A1) and the inequality $ \mathbb{E}(\xi_{i}^2 | \mathcal{F}_{i-1})  \leq   12\, \epsilon^{2}$ (cf. Lemma \ref{l11}),  we have, for any constant $\alpha \in (0,1)$ and all $0 \leq \lambda \leq \alpha\, \epsilon^{-1} ,$
\[
\left|\Psi _n(\lambda ) - \frac{\lambda ^2}2\left\langle X\right\rangle _n \right|  \leq  c_\alpha \lambda^3 \epsilon \langle X\rangle_n.
\]
Combining this inequality with condition (A2), we get, for any constant $\alpha \in (0,1)$ and all $0 \leq \lambda \leq \alpha\, \epsilon^{-1} ,$
\[
\left|\Psi _n(\lambda ) - \frac{\lambda ^2}{2} \right| \leq
2 \,c_\alpha  \lambda^3 \epsilon + \frac{\lambda^2 \delta^2}2,
\]
which completes the proof of Lemma \ref{LEMMA-1-2}.\hfill\qed

\section{\textbf{Proof of Theorem  \ref{th0}}\label{sec3.1}}
For $0\leq x < 1$, the assertion follows from Theorem \ref{th4}.
It remains to prove Theorem \ref{th0} for any $\alpha \in (0,1)$ and all $1 \leq x \leq \alpha\, \epsilon^{-1}$.
Changing the probability measure according to (\ref{f21}),  we have, for all $0\leq \lambda <\, \epsilon^{-1},$
\begin{eqnarray}
\mathbb{P}(X_n>x) &= & \mathbb{E}_\lambda \left( Z_n (\lambda)^{-1}\mathbf{1}_{\{ X_n>x \}} \right) \nonumber\\
&= & \mathbb{E}_\lambda \left(\exp \left\{-\lambda X_n+\Psi _n(\lambda )\right\} \mathbf{1}_{\{ X_n>x \}} \right) \nonumber\\
&= & \mathbb{E}_\lambda \left(\exp \left\{
-\lambda Y_n(\lambda)-\lambda B_{n}(\lambda)+\Psi _n(\lambda)\right\} \mathbf{1}_{\{Y_n(\lambda)+B_{n}(\lambda)>x\}} \right).\label{f32}
\end{eqnarray}
Let
$\overline{\lambda}=\overline{\lambda}(x)$ be the largest
solution of the equation
\begin{equation}\label{f33}
\lambda +\lambda \delta^{2} +c_\alpha \lambda^{2}\epsilon  =x,
\end{equation}
where $c_\alpha$ is given by inequality (\ref{f25sa}).
The definition of $\overline{\lambda}$ implies that there exist $c_{\alpha,0}, c_{\alpha,1} >0$ such that, for all $1 \leq x \leq \alpha\, \epsilon^{-1} ,$
\begin{equation}\label{f34}
c_{\alpha,0}\,x \leq \overline{\lambda}=\frac{2x}{\sqrt{(1+\delta^2)^2+ 4 c_\alpha x\epsilon}+1+\delta^2} \leq  x
\end{equation}
and
\begin{equation}\label{f35}
\overline{\lambda}=x - c_{\alpha,1}|\theta|(x^2\epsilon +x \delta ^2) \in [c_{\alpha,0}, \alpha\, \epsilon^{-1}\,].
\end{equation}
 From (\ref{f32}), using Lemmas \ref{LEMMA-1-1}, \ref{LEMMA-1-2} and equality (\ref{f33}),
 we obtain, for all $1 \leq x \leq \alpha\, \epsilon^{-1} ,$
\begin{equation}\label{gfdgf}
\mathbb{P}(X_n>x)\leq  e^{ c_{\alpha,2} \,(\overline{\lambda}^3\epsilon+\overline{\lambda}^2\delta^2) -\overline{\lambda}^2/2}\mathbb{E}_{\overline{\lambda}}\left(e^{-%
\overline{\lambda}Y_n(\overline{\lambda})}\mathbf{1}_{\{ Y_n(\overline{\lambda})>0\}} \right).
\end{equation}
It is easy to see that
\begin{equation}\label{f37d1}
\mathbb{E}_{\overline{\lambda}} \left( e^{-%
\overline{\lambda}Y_n(\overline{\lambda})}\mathbf{1}_{\{ Y_n(\overline{\lambda})>0\}}\right)= \int_{0}^{\infty} \overline{\lambda} e^{-\overline{\lambda} y}  \mathbb{P}_{\overline{\lambda}}(0 < Y_n(\overline{\lambda})\leq y  ) dy.
\end{equation}
Similarly, for a standard gaussian random variable $\mathcal{N}$, we have
\begin{equation}\label{f37d2}
\mathbb{E} \left( e^{-%
\overline{\lambda}\mathcal{N}}\mathbf{1}_{\{ \mathcal{N}>0\}} \right)= \int_{0}^{\infty} \overline{\lambda} e^{-\overline{\lambda} y}   \mathbb{P} (0 < \mathcal{N} \leq y  ) dy.
\end{equation}
From (\ref{f37d1}) and (\ref{f37d2}), it follows
\begin{eqnarray}
\left|\mathbb{E}_{\overline{\lambda}}\left(e^{-%
\overline{\lambda}Y_n(\overline{\lambda})}\mathbf{1}_{\{ Y_n(\overline{\lambda})>0\}}\right)- \mathbb{E} \left( e^{-%
\overline{\lambda}\mathcal{N}}\mathbf{1}_{\{ \mathcal{N}>0\}}\right) \right| \leq 2\sup_y \bigg| \mathbb{P}_{\overline{\lambda}} (Y_n(\overline{\lambda} )\leq y)-\Phi (y) \bigg|.\nonumber
\end{eqnarray}
Using Lemma \ref{LEMMA4}, we obtain the following bound: for all $
1 \leq  x  \leq \alpha\,  \epsilon ^{-1},$
\begin{equation}\label{f38}
\left|\mathbb{E}_{\overline{\lambda}}\left(e^{-%
\overline{\lambda}Y_n(\overline{\lambda})}\mathbf{1}_{\{ Y_n(\overline{\lambda})>0\}}\right)- \mathbb{E}\left( e^{-%
\overline{\lambda}\mathcal{N}}\mathbf{1}_{\{ \mathcal{N}>0\}}\right)\right| \leq  c\left( \overline{\lambda}\epsilon +\epsilon \left| \log \epsilon
 \right| +\delta \right).
\end{equation}
From (\ref{gfdgf}) and (\ref{f38}) we find that, for all $
1 \leq  x  \leq \alpha \,  \epsilon ^{-1}, $
\[
\mathbb{P}(X_n>x)\leq e^{ c_{\alpha,2} \,(\overline{\lambda}^3\epsilon+\overline{\lambda}^2\delta^2) -\overline{\lambda}^2/2}\Bigg( \mathbb{E}  \Big( e^{- \overline{\lambda}\mathcal{N}}\mathbf{1}_{\{ \mathcal{N}>0\}}  \Big)+c\left( \overline{\lambda}\epsilon +\epsilon \left| \log   \epsilon
 \right| +\delta \right)  \Bigg).
\]
Since
\begin{equation}\label{fsphi}
e^{- \lambda^2/2}\mathbb{E} \left( e^{-
 \lambda\mathcal{N}}\mathbf{1}_{\{ \mathcal{N}>0\}} \right) =\frac{1}{\sqrt{2\pi}}\int_0^{\infty}e^{-(y+\lambda)^2/2}
  dy=1-\Phi \left(  \lambda\right)
\end{equation}
and, for all $\lambda\geq c_{\alpha,0},$
\begin{equation}
1-\Phi \left(  \lambda\right)
\geq
\frac 1{\sqrt{2 \pi}(1+ \lambda)}\ e^{- \lambda^2/2}
\geq
\frac{c_{\alpha,0}}{\sqrt{2 \pi}(1+ c_{\alpha,0})} \frac 1{ \lambda }e^{- \lambda^2/2}
\label{f39}
\end{equation}
(see Feller \cite{F71}), we obtain the following upper bound on tail probabilities:
for all $1\leq x \leq \alpha  \, \epsilon ^{-1}  ,$
\begin{eqnarray}
\frac{\mathbb{P}(X_n>x)}{1-\Phi \left( \overline{\lambda}\right)}
&\leq& e^{ c_{\alpha,2} \,(\overline{\lambda}^3\epsilon+\overline{\lambda}^2\delta^2) } \left(\, 1+ c_{\alpha,3}\, (\, \overline{\lambda}^2 \epsilon +\overline{\lambda}\epsilon \left| \log \epsilon \right| +\overline{\lambda}\delta \, ) \right).  \label{f40}
\end{eqnarray}
Next, we would like to compare $1-\Phi (\overline{\lambda})$ with $1-\Phi (x)$.
  By (\ref{f34}), (\ref{f35}) and (\ref{f39}),  we get
\begin{eqnarray}
   1 \leq \frac{\int_{\overline{\lambda}}^{ \infty}\exp\{- t^2/2 \}d t}{\int_{x}^{ \infty}\exp\{- t^2/2 \}d t}   &= &
  1+\frac{\int_{\overline{\lambda}}^{x}\exp\{ -t^2/2 \} d t}{\int_{x}^{ \infty}\exp\{-t^2/2\}d t}\nonumber\\
   & \leq & 1+c_{\alpha,4}x(x-\overline{\lambda}) \exp\{ (x^2-\overline{\lambda}^2)/2 \}\nonumber\\
   & \leq & \exp\{ c_{\alpha,5}\, (x^3 \epsilon + x^2 \delta^2)\}.\label{f41}
\end{eqnarray}
So, we find that
\begin{equation}\label{f42}
1-\Phi \left( \overline{\lambda}\right) =\Big( 1-\Phi (x)\Big)\exp \left\{ |\theta_{1}| c_{\alpha,5}\, (
x^3 \epsilon+ x^2 \delta^2 ) \right\}.
\end{equation}
Implementing (\ref{f42}) in (\ref{f40}) and using (\ref
{f34}), we obtain, for all $1\leq x \leq \alpha \, \epsilon ^{-1}, $
\begin{eqnarray*}
\frac{\mathbb{P}(X_n>x)}{1-\Phi \left( x\right) }&\leq& \exp\{ c_{\alpha,6} ( x^3 \epsilon + x^2 \delta^2 ) \}
\left(\frac{}{} 1+c_{\alpha,7}\left( x^2 \epsilon +x\epsilon \left| \log  \epsilon
 \right| +x\delta \right) \right)\\
&\leq & \exp\{ c_{\alpha,6} ( x^3 \epsilon + x^2 \delta^2 )  \}\left(\frac{}{} 1+c_{\alpha,7} \, x^2\epsilon \right)
\left(\frac{}{} 1+c_{\alpha,7}\, x \left( \epsilon \left| \log \epsilon
 \right| + \delta \right) \right) \nonumber\\
&\leq & \exp\{ c_{\alpha,8} ( x^3 \epsilon + x^2 \delta^2 )  \}
\left(\frac{}{} 1+c_{\alpha,7} \, x\left( \epsilon \left| \log  \epsilon
 \right| + \delta \right) \right).
\end{eqnarray*}
Taking $c_{\alpha}=\max\{c_{\alpha,7}, c_{\alpha,8}\}$, we prove the first assertion of Theorem \ref{th0}.
The second assertion follows from the first one applied to the martingale  $(-X_k)_{k=0,...,n}$.

\section{\textbf{Proof of Theorem  \ref{th1}}\label{sec3.2}}
For $0\leq x < 1$, the assertion follows from Theorem \ref{th4}.
It remains to prove Theorem \ref{th1} for $1 \leq x \leq \alpha_0 \epsilon^{-1}$, where $\alpha_0 >0$ is  an absolute constant.
 Let
$\underline{\lambda}=\underline{\lambda}(x)$ be the smallest solution of the equation
\begin{equation}\label{f44}
\lambda -\lambda \delta^{2} -c_{1/2} \lambda^{2}\epsilon   =x,
\end{equation}
where $c_\alpha$ is given by inequality (\ref{f25sa}).
The definition of $\underline{\lambda}$ implies that, for all $1 \leq x \leq 0.01 c_{1/2}^{-1} \epsilon^{-1}, $ it holds
\begin{equation}\label{f45}
x\leq \underline{\lambda}=\frac{2x}{1-\delta^2 + \sqrt{(1-\delta^2)^2-4c_{1/2} x\epsilon }} \leq 2 \,x
\end{equation}
and
\begin{equation}\label{f46}
\underline{\lambda}=x+c_{0}|\theta|(x^2\epsilon +x \delta ^2) \in [1, 0.02\, c_{1/2}^{-1}  \epsilon^{-1}].
\end{equation}
 From (\ref{f32}), using Lemmas \ref{LEMMA-1-1},  \ref{LEMMA-1-2} and equality (\ref{f44}),
  we obtain, for all $1 \leq x \leq 0.01 c_{1/2}^{-1} \epsilon^{-1}, $
\begin{equation}\label{jknjssa}
\mathbb{P}(X_n>x)\geq e^{ - c_{1} \,(\underline{\lambda}^3\epsilon+\underline{\lambda}^2\delta^2) -\underline{\lambda}^2/2} \mathbb{E}_{\underline{\lambda}} \left(e^{-%
\underline{\lambda}Y_n(\underline{\lambda})}\mathbf{1}_{\{ Y_n(\underline{\lambda})>0 \}} \right).
\end{equation}

In the subsequent we distinguish two cases. First, let $1\leq \underline{\lambda}\leq \alpha_1 \min\{ \epsilon ^{-1/2} , \delta^{-1} \}$, where $\alpha_1>0$ is a small absolute constant
whose value will be given later.
Note that inequality (\ref{f38}) can be established with $\overline{\lambda}$ replaced by $\underline{\lambda}$, which, in turn,
implies
\[
\mathbb{P}(X_n>x)\geq e^{ - c_{1} \,(\underline{\lambda}^3\epsilon+\underline{\lambda}^2\delta^2) -\underline{\lambda}^2/2} \bigg( \mathbb{E} \left(e^{-\underline{\lambda}  \mathcal{N}}\mathbf{1}_{\{ \mathcal{N}>0\}} \right)-c_2\left( \underline{\lambda}\epsilon +\epsilon \left| \log   \epsilon
 \right| +\delta \right)  \bigg).
\]
By (\ref{fsphi}) and (\ref{f39}),
we obtain the following lower bound on tail probabilities:
\begin{eqnarray}
\frac{\mathbb{P}(X_n>x)}{1-\Phi \left( \underline{\lambda}\right)}
&\geq& e^{ - c_{1} \,(\underline{\lambda}^3\epsilon+\underline{\lambda}^2\delta^2) } \left(\frac{}{} 1-c_2\left( \underline{\lambda}^2 \epsilon +\underline{\lambda}\epsilon \left| \log  \epsilon
 \right| +\underline{\lambda}\delta \right) \right).  \label{f51}
\end{eqnarray}
Taking  $\alpha_1 =  (8c_{2})^{-1}$,
we deduce that, for all $1\leq \underline{\lambda} \leq \alpha_1 \min \{\epsilon^{-1/2}, \delta^{-1} \} $,
\begin{eqnarray}\label{f55f}
1-c_2\left( \underline{\lambda}^2 \epsilon +\underline{\lambda}\epsilon \left| \log  \epsilon
 \right| +\underline{\lambda}\delta \right)&\geq &  \exp\left\{-2c_2\left( \underline{\lambda}^2 \epsilon +\underline{\lambda}\epsilon \left| \log  \epsilon
 \right| +\underline{\lambda}\delta \right) \right\} .
\end{eqnarray}
Implementing (\ref{f55f}) in (\ref{f51}), we obtain
\begin{eqnarray}
\frac{\mathbb{P}(X_n>x)}{1-\Phi \left( \underline{\lambda}\right)}&\geq& \exp \bigg \{ -c_3 \left( \underline{\lambda}^3 \epsilon +\underline{\lambda }\epsilon \left| \log  \epsilon
 \right| +\underline{\lambda}\delta +\underline{\lambda}^2\delta^2 \right)\bigg\}   \label{f54}
\end{eqnarray}
which is valid for all $1\leq \underline{\lambda} \leq \alpha_1 \min \{\epsilon^{-1/2}, \delta^{-1} \} $.

Next, we consider the case of $\alpha_1 \min\{ \epsilon ^{-1/2} , \delta^{-1} \} \leq \underline{\lambda} \leq \alpha_0 \epsilon ^{-1}$ and $\delta \leq \alpha_0$. Let $K \geq 1$ be an absolute constant, whose exact value will be chosen later.
It is easy to see that
\begin{eqnarray}\label{jknjsta}
\mathbb{E}_{\underline{\lambda}} \left(e^{-\underline{\lambda}Y_n(\underline{\lambda})}\mathbf{1}_{\{ Y_n(\underline{\lambda})>0 \}} \right) &\geq& \mathbb{E}_{\underline{\lambda}} \Big(e^{-\underline{\lambda}Y_n(\underline{\lambda})}\mathbf{1}_{\{0< Y_n(\underline{\lambda})\leq  K \gamma \}} \Big) \nonumber\\
 &\geq&e^{-\underline{\lambda} K \gamma}\mathbb{P}_{\underline{\lambda}} \Big(0< Y_n(\underline{\lambda})\leq  K \gamma \Big),
\end{eqnarray}
where $\gamma = \underline{\lambda}\epsilon +  \epsilon |\log \epsilon|  +\delta \le 4 \alpha_0^{1/2},$ if $\alpha_0\leq 1.$
From Lemma \ref{LEMMA4}, we have
\begin{eqnarray*}
\mathbb{P}_{\underline{\lambda}} \Big(0< Y_n(\underline{\lambda})\leq  K \gamma \Big) &\geq&   \mathbb{P}  \Big( 0<  \mathcal{N}
\leq  K \gamma  \Big)  - c_5 \gamma  \\
  &\geq& K \gamma  e^{-K^2  \gamma ^2/2}   - c_5 \gamma\\
  &\geq& \left( K e^{-8 K^2  \alpha_0}   - c_5 \right) \gamma.
\end{eqnarray*}
Taking $\alpha_0=1/(16 K^2)$, we find that $$\mathbb{P}_{\underline{\lambda }} \Big(0< Y_n(\underline{\lambda})\leq  K \gamma \Big)  \geq
    \left( \frac12 K   - c_5 \right) \gamma.$$
Letting $K\geq   8c_5 $, it follows that
$$\mathbb{P}_{\underline{\lambda}} \Big(0< Y_n(\underline{\lambda })\leq  K \gamma \Big)  \geq \frac38 K\gamma \geq \frac38 K  \frac{  \max \left\{ \underline{\lambda }^2\epsilon  , \underline{\lambda } \delta  \right\}} { \underline{\lambda } }.$$
Choosing $K=\max \Big\{8 c_5,  \frac{8\alpha_1^{-2}}{3\sqrt{\pi}} \Big\}$ and taking into account that
$ \alpha_1 \min\{ \epsilon ^{-1/2},  \delta^{-1} \} \leq \underline{\lambda} \leq \alpha_0 \epsilon^{-1}$,
we deduce that
\begin{eqnarray*}
\mathbb{P}_{\underline{\lambda }} \Big(0< Y_n(\underline{\lambda }) \leq  K \gamma \Big)  \geq    \frac{1}{\sqrt{\pi}\underline{\lambda} } .
\label{jknjstb}
\end{eqnarray*}
Since the inequality
$\frac{1}{\sqrt{\pi}\lambda }  e^{- \lambda^2/2}  \geq 1-\Phi \left(  \lambda \right)$ is valid for $\lambda \geq 1$ (see Feller \cite{F71}),
it follows that, for all $ \alpha_1 \min\{ \epsilon ^{-1/2},  \delta^{-1} \} \leq \underline{\lambda} \leq \alpha_0 \epsilon^{-1}$,
\begin{eqnarray}\label{dfac}
\mathbb{P}_{\underline{\lambda}} \Big(0< Y_n(\underline{\lambda})\leq  K \gamma \Big)   \geq \bigg( 1-\Phi \left(  \underline{\lambda} \right) \bigg)e^{\underline{\lambda}^2/2} .
\end{eqnarray}
From (\ref{jknjssa}), (\ref{jknjsta}) and (\ref{dfac}), we obtain
 \begin{eqnarray}
\frac{\mathbb{P}(X_n>x)}{1-\Phi \left( \underline{\lambda}\right)}&\geq& \exp \bigg \{ -c_{\alpha_0, 6} \left( \underline{\lambda}^3 \epsilon   + \underline{\lambda} \epsilon |\log \epsilon|+\underline{\lambda}\delta +\underline{\lambda}^2\delta^2 \right)\bigg\}  \label{fgj53}
\end{eqnarray}
which is valid for all $ \alpha_1 \min\{ \epsilon ^{-1/2},  \delta^{-1} \} \leq \underline{\lambda} \leq \alpha_0 \epsilon^{-1}$.

Putting (\ref{f54}) and (\ref{fgj53}) together, we obtain, for all $1\leq \underline{\lambda} \leq \alpha_0 \epsilon^{-1}$ and $\delta \leq \alpha_0$,
\begin{eqnarray}\label{ft52}
\frac{\mathbb{P}(X_n>x)}{1-\Phi \left( \underline{\lambda}\right)}&\geq& \exp \bigg \{ -c_{\alpha_0, 7} \left( \underline{\lambda}^3 \epsilon  +  \underline{\lambda}\epsilon \left| \log  \epsilon  \right|+\underline{\lambda} \delta + \underline{\lambda}^2\delta^2\right)\bigg\}  .
\end{eqnarray}
As in the proof of Theorem \ref{th0}, we now compare $1-\Phi (\underline{\lambda})$ with  $1-\Phi (x)$. By a similar argument as in (\ref{f41}), we have
\begin{equation}\label{f53}
1-\Phi \left( \underline{\lambda}\right) =\Big( 1-\Phi (x)\Big)  \exp \left\{ - |\theta| c_3\, (
x^3 \epsilon+ x^2 \delta^2 ) \right\} .
\end{equation}
Combining (\ref{f45}), (\ref{ft52}) and (\ref{f53}), we obtain, for all $1\leq x \leq \alpha_0 \epsilon^{-1}$ and $\delta \leq \alpha_0$,
 \begin{equation}
\frac{\mathbb{P}(X_n>x)}{1-\Phi \left( x\right)}  \geq  \exp \bigg \{ -c_{\alpha_0, 8} \left( x^3 \epsilon  + x\epsilon \left| \log  \epsilon  \right|+x \delta + x^2\delta^2\right)\bigg\}
\end{equation}
 which gives the first conclusion of Theorem \ref{th1}.
The second conclusion follows from the first one applied to the martingale  $(-X_k)_{k=0,...,n}$.\hfill\qed

\section{\textbf{Proof of Lemma \ref{LEMMA4}}\label{sec4}}
The proof of Lemma \ref{LEMMA4} is a refinement of Lemma 3.3 of Grama and  Haeusler \cite{GH00} where it is assumed that $|\eta_{i}|\leq 2\epsilon$,  which is a  particular case of condition  (A1). Compared to the case where $\eta_{i}$ are bounded, the main challenge of our proof comes from the control of $I_1$ defined in (\ref{D-1}) below.

In this section, $\alpha$ denotes a positive absolute number satisfying $\alpha \in (0,1)$,    $\vartheta$ denotes a real number satisfying $0\leq \vartheta \leq 1 ,$ which is different from  $\theta$, and
$\varphi(t)$ denotes the density function of the standard normal distribution.  For the sake of simplicity, we also denote $Y(\lambda), Y_{n}(\lambda)$ and $ \eta(\lambda)$ by $Y, Y_{n}$ and $\eta$, respectively.
We want to obtain a rate of convergence in the central limit
theorem for the conjugate martingale $Y=(Y_k,\mathcal{F}_k)_{k=1,...,n},$
where $Y_k=\sum_{i=1}^k\eta _i .$
Denote the quadratic characteristic of the conjugate
martingale $Y$ by $\left\langle Y \right\rangle_{k}= \sum_{i\leq k}\mathbb{E}_{\lambda}(\eta_{i}^2 | \mathcal{F}_{i-1})$, and set
$\Delta \left\langle Y \right\rangle_{k}=\mathbb{E}_{\lambda}(\eta_{k}^2 | \mathcal{F}_{k-1}).$ It is easy to see that, for $k=1,...,n$,
\begin{eqnarray}
\Delta \left\langle Y \right\rangle_{k} &=&   \mathbb{E}_{\lambda}\left((\xi_{k}-b_{k}(\lambda))^2  | \mathcal{F}_{k-1} \right) \nonumber\\
&= & \frac{\mathbb{E}(\xi
_k^2e^{\lambda \xi _k}|\mathcal{F}_{k-1})}{\mathbb{E}(e^{\lambda \xi _k}|\mathcal{F}%
_{k-1})}-\frac{\mathbb{E}(\xi _ke^{\lambda \xi _k}|\mathcal{F}_{k-1})^2}{\mathbb{E}(e^{\lambda
\xi _k}|\mathcal{F}_{k-1})^2}. \label{f24}
\end{eqnarray}

Since $\mathbb{E}(e^{\lambda\xi_{i}} | \mathcal{F}_{i-1} ) \geq 1 $ and $|\eta_{i}|^k \leq 2^{k-1}(|\xi_{i}|^k+ \mathbb{E}_{\lambda}(|\xi_{i}||\mathcal{F}_{i-1})^k)$, using condition (A1) and Lemma \ref{l11}, we obtain, for all $k\geq 3$ and all $0\leq \lambda \leq \frac14\, \epsilon^{-1}$,
\begin{eqnarray*}
\mathbb{E}_{\lambda}\left(|\eta_{i}|^{k} | \mathcal{F}_{i-1} \right) &\leq &2^{k-1} \mathbb{E}_{\lambda}\left(|\xi_{i}|^k + \mathbb{E}_{\lambda}(|\xi_{i}| | \mathcal{F}_{i-1})^k | \mathcal{F}_{i-1} \right)\\
&\leq & 2^k \mathbb{E}_{\lambda}\left(|\xi_{i}|^k  | \mathcal{F}_{i-1} \right) \\
&\leq & 2^k \mathbb{E} \left(|\xi_{i}|^k \exp\{|\lambda\xi_{i}|\} | \mathcal{F}_{i-1} \right) \\
&\leq & c 2^{k}k!\epsilon^{k-2}\mathbb{E} \left(\xi_{i}^2 | \mathcal{F}_{i-1} \right)\,.
\end{eqnarray*}
Using Taylor's expansion for $e^x$ and Lemma 1, we have, for all $0 \leq \lambda \leq \frac14\, \epsilon^{-1} ,$
\begin{eqnarray}\label{f56}
\left| \Delta \left\langle Y\right\rangle _k-\Delta \left\langle
X\right\rangle _k\right| &\leq &\left| \frac{\mathbb{E}(\xi _k^2e^{\lambda \xi _k}|%
\mathcal{F}_{k-1})}{\mathbb{E}(e^{\lambda \xi _k}|\mathcal{F}_{k-1})}-\mathbb{E}(\xi _k^2|%
\mathcal{F}_{k-1})\right| +\left| \frac{\mathbb{E}(\xi _k e^{\lambda \xi _k}|\mathcal{F%
}_{k-1})^2}{\mathbb{E}(e^{\lambda \xi _k}|\mathcal{F}_{k-1})^2}\right|  \nonumber \\
&\leq &\left|\mathbb{E}(\xi _k^2e^{\lambda \xi _k}|\mathcal{F}_{k-1})-\mathbb{E}(\xi _k^2|\mathcal{F}_{k-1})\mathbb{E}(e^{\lambda \xi _k}|\mathcal{F}_{k-1})\right| \nonumber \\
&& + \mathbb{E}(\xi_{k} e^{\lambda \xi _k}|\mathcal{F}_{k-1})^2 \nonumber\\
&\leq & \sum_{l=1}^{\infty}|\mathbb{E}(\xi_{k}^{l+2} | \mathcal{F}_{k-1})|\frac{\lambda^l}{l !} + \Delta \langle X\rangle_{k}
\sum_{l=1}^{\infty}|\mathbb{E}(\xi_{k}^l | \mathcal{F}_{k-1})|\frac{\lambda^l}{l !} \nonumber\\
&& +\left(\sum_{l=1}^{\infty}|\mathbb{E}(\xi_{k}^{l+1} | \mathcal{F}_{k-1})|\frac{ \lambda ^l}{l !}\right)^2 \nonumber\\
&\leq &  c \lambda\epsilon\, \Delta \langle X\rangle_{k} .
\end{eqnarray}
Therefore,
\[
|\langle Y\rangle_{n} - 1| \leq |\langle Y\rangle_{n} - \langle X\rangle_{n}| + |\langle X\rangle_{n}-1| \leq c \lambda\epsilon \langle X\rangle_{n}  + \delta^2 .
\]
Thus the martingale $Y$
satisfies the following conditions (analogous to conditions (A1) and (A2)): for all $0 \leq \lambda \leq \frac14 \, \epsilon^{-1}$,
\begin{description}
\item[(B1)] $\mathbb{E}_{\lambda}(|\eta_{i}|^{k} | \mathcal{F}_{i-1} )\leq c_{k} \epsilon^{k-2}\mathbb{E}(\xi_{i}^2 | \mathcal{F}_{i-1} ) , \ \ \ 5 \geq k \geq 3;$
\item[(B2)] $\; \left| \left\langle Y\right\rangle _n-1\right| \leq  c( \lambda\epsilon + \delta^2 )  .$
\end{description}

We first prove Lemma \ref{LEMMA4} for $1\leq \lambda < \epsilon ^{-1}$. Without loss of generality,
we can assume that  $1\leq \lambda \leq  \frac14 \, \epsilon ^{-1}$, otherwise we take $c \geq 4$ in the assertion of the lemma.
Set $T=1+\delta ^2$ and introduce a modification of the
quadratic characteristic $\left\langle X\right\rangle $ as follows:
\begin{equation}
V_k=\left\langle X\right\rangle _k\mathbf{1}_{\{k<n\}}+T\mathbf{1}_{\{k=n\}}.
\label{RA-4}
\end{equation}
Note that $V_0=0,$ $V_n=T$ and that $(V_k,\mathcal{F}_k)_{k=0,...,n}$ is a
predictable process. Set $\gamma =\lambda \epsilon +\delta,$ where $\lambda \in [1, \epsilon^{-1})$. Let $%
c_{*}\geq 4$ be a ``free'' absolute constant, whose exact value will be
chosen later. Consider the non-increasing discrete time predictable process
$A_k=c_{*}^2\gamma ^2+T-V_k, k=1,...,n.$ For any fixed $u\in \mathbb{R}$ and any $x\in \mathbb{R} $ and $
y > 0,$ set for brevity,
\begin{equation}
\Phi _u(x,y)=\Phi \left( (u-x)/\sqrt{y}\right) .  \label{RA-7}
\end{equation}

In the proof we make use of the following two assertions, which
can be found in Bolthausen's paper \cite{Bo82}.
\begin{lemma}\cite{Bo82}
\label{LEMMA-APX-1}Let $X$ and $Y$ be random variables. Then
\[
\sup_u\left| \mathbb{P}\left( X\leq u\right) -\Phi \left( u \right) \right|
\leq c_1\sup_u\left| \mathbb{P}\left( X+Y\leq u\right) -\Phi \left( u \right)
\right| +c_2\left\| \mathbb{E}\left( Y^2|X\right) \right\| _\infty ^{1/2}.
\]
\end{lemma}
\begin{lemma} \cite{Bo82}
\label{LEMMA-APX-2}Let $G(x)$ be an integrable function of bounded variation,
$X$ be a random variable and $a,$ $b>0$ are real numbers. Then
\[
\mathbb{E} \, G\left( \frac{X+a}b\right) \leq c_1\sup_u\left| \mathbb{P}\left( X\leq u\right)
-\Phi \left( u\right) \right| +c_2\, b.
\]
\end{lemma}

Let $\mathcal{N}_{c_{*}^2\gamma ^2}=\mathcal{N}(0,c_{*}\gamma)$ be a normal random variable  independent of $Y_n$.
Using a well-known smoothing procedure (which employs Lemma \ref{LEMMA-APX-1}),
we get
\begin{eqnarray}
 \sup_u\left| \mathbb{P}_{\lambda}(Y_n\leq u)  -  \Phi (u)\right| &\leq& c_1\sup_u\left|
\mathbb{E}_{\lambda}\Phi _u(Y_n,A_n)- \Phi (u)\right|
+c_2\gamma \nonumber\\
&\leq& c_1\sup_u\left| \mathbb{E}_{\lambda}\Phi _u(Y_n,A_n)-\mathbb{E}_{\lambda}\Phi _u(Y_0,A_0)\right| \nonumber\\
& & +\, c_1\sup_u\left| \mathbb{E}_{\lambda}\Phi _u(Y_0,A_0)-\Phi (u)\right| +c_2\gamma   \nonumber\\
&=& c_1\sup_u\left| \mathbb{E}_{\lambda}\Phi _u(Y_n,A_n)-\mathbb{E}_{\lambda}\Phi _u(Y_0,A_0)\right| \nonumber\\
& & +\, c_1\sup_u\left| \Phi \left(\frac{u}{\sqrt{c_{*}^2\gamma^2+T}}  \right)-\Phi (u)\right| +c_2\gamma    \nonumber\\
&\leq& c_1\sup_u\left| \mathbb{E}_{\lambda}\Phi _u(Y_n,A_n)-\mathbb{E}_{\lambda}\Phi _u(Y_0,A_0)\right|  +c_3 \gamma, \label{llaa}
\end{eqnarray}
where
\[
\mathbb{E}_{\lambda}\Phi _u(Y_n,A_n)=\mathbb{P}_{\lambda}(Y_n+\mathcal{N}_{c_{*}^2\gamma ^2}\leq u) \ \ \mbox{and} \ \
\mathbb{E}_{\lambda}\Phi _u(Y_0,A_0)=\mathbb{P}_{\lambda}(\mathcal{N}_{c_{*}^2\gamma^2+T}\leq u).
\]
By simple telescoping, we find that
\[
\mathbb{E}_\lambda \Phi _u(Y_n,A_n)-\mathbb{E}_\lambda \Phi _u(Y_0,A_0)=\mathbb{E}_\lambda
\sum_{k=1}^n\bigg( \Phi _u(Y_k,A_k)-\Phi _u(Y_{k-1},A_{k-1})\bigg) .
\]
From this, taking into account that $(\eta _i,\mathcal{F}_i)_{i=0,...,n}$ is
a $\mathbb{P}_\lambda $-martingale and that
\[
\frac{\partial ^2}{\partial x^2}\Phi _u(x,y)=2\frac \partial {\partial
y}\Phi _u(x,y),
\]
we obtain
\begin{equation}
\mathbb{E}_\lambda \Phi _u(Y_n,A_n)-\mathbb{E}_\lambda \Phi _u(Y_0,A_0)=I_1+I_2-I_3,
\label{lla}
\end{equation}
where
\begin{eqnarray}
I_1&=&\mathbb{E}_\lambda \sum_{k=1}^n\Bigg( \frac{}{} \Phi _u(Y_k,A_k)-\Phi
_u(Y_{k-1},A_k)  \nonumber \\
&& \ \ \ \ \ \ \ \ \ \ \ \ \   -\frac \partial {\partial x}\Phi
_u(Y_{k-1},A_k)\eta _k-\frac 12\frac{\partial ^2}{\partial x^2}\Phi
_u(Y_{k-1},A_k)\eta _k^2\Bigg) ,  \label{D-1}\\
I_2 &=& \frac 12\mathbb{E}_\lambda \sum_{k=1}^n\frac{\partial ^2}{\partial x^2}\Phi
_u(Y_{k-1},A_k)\bigg(\Delta \left\langle Y\right\rangle _k-\Delta V_k \bigg),\quad \quad \
\label{D-2}
\end{eqnarray}
\begin{equation}
I_3=\mathbb{E}_\lambda \sum_{k=1}^n\left( \Phi _u(Y_{k-1},A_{k-1})-\Phi
_u(Y_{k-1},A_k)-\frac \partial {\partial y}\Phi _u(Y_{k-1},A_k)\Delta
V_k\right) .  \label{D-3}
\end{equation}
We now give estimates of $I_1,$ $I_2$ and $I_3.$ To shorten notations, set
\[T_{k-1}= (u-Y_{k-1})/\sqrt{A_k}.\]

\emph{\textbf{a)} Control of }$I_1 .$ \quad Using a three-term Taylor's expansion, we have
\begin{equation}
I_1=-\mathbb{E}_\lambda \sum_{k=1}^n\frac 1{6A_k^{3/2}}\varphi ^{\prime \prime }\left(T_{k-1}-\frac{\vartheta_{k}\eta_{k}}{ \sqrt{A_{k}}} \right) \eta _k^3\ .  \label{CC-1}
\end{equation}
In order to bound $\varphi ^{\prime \prime } (\cdot)$ we distinguish two cases as follows.

\emph{Case 1}: $|\eta_{k}/\sqrt{A_{k}}| \leq |T_{k-1}|/2.$ In this case, by the inequality $\varphi''(t)\leq\varphi(t)(1+t^2)$, it follows
\begin{eqnarray*}
  \left|\varphi ^{\prime \prime }\left(T_{k-1}-\frac{\vartheta_{k}\eta_{k}}{ \sqrt{A_{k}}} \right) \right| &\leq & \varphi\left(T_{k-1}-\frac{\vartheta_{k}\eta_{k}}{ \sqrt{A_{k}}} \right)\left(1+\left(T_{k-1}-\frac{\vartheta_{k}\eta_{k}}{ \sqrt{A_{k}}} \right)^2 \right) \\
  &\leq & \sup_{|t-T_{k-1}|\leq |T_{k-1}|/2}\varphi(t)(1+t^2) \\
 & \leq & \varphi(T_{k-1}/2)(1+4T_{k-1}^2).
\end{eqnarray*}
Define $g_{1}(t)=\sup_{|t-z|\leq 3}f_{1}(z),$ where $f_{1}(t) = \varphi(t/2)(1+4t^2).$
It is easy to see that $g_{1}(t)$ is a symmetric integrable function of bounded variation, non-increasing in $t\geq 0$. Therefore,
\begin{eqnarray}\label{fklm59}
\left|\varphi ^{\prime \prime }\left(T_{k-1}-\frac{\vartheta_{k}\eta_{k}}{ \sqrt{A_{k}}} \right) \right| \mathbf{1}_{\left\{|\eta_{k}/\sqrt{A_{k}}| \leq |T_{k-1}|/2 \right\}} \leq  g_{1}(T_{k-1}) .
\end{eqnarray}

\emph{Case 2}: $|\eta_{k}/\sqrt{A_{k}}| > |T_{k-1}|/2$. Since $|\varphi''(t)| \leq 2$, it follows that
\begin{eqnarray}\label{fktm59}
\left|\varphi ^{\prime \prime }\left(T_{k-1}-\frac{\vartheta_{k}\eta_{k}}{ \sqrt{A_{k}}} \right) \right|\mathbf{1}_{\left\{|\eta_{k}/\sqrt{A_{k}}| > |T_{k-1}|/2 \right\}} \leq  2 \left( \mathbf{1}_{\{|T_{k-1}|<2\}}+ \frac{4\eta_{k}^2}{T_{k-1}^2 A_{k}} \mathbf{1}_{\{|T_{k-1}|\geq2\}} \right).
\end{eqnarray}

Now we bound the conditional expectation of $|\eta_{k}|^k$. Using condition (B1), we have
\[
\mathbb{E}_{\lambda}(|\eta_{k}|^3 | \mathcal{F}_{k-1}) \leq c\, \epsilon \Delta \langle X \rangle_{k} \quad
\textrm{and}\quad
\mathbb{E}_{\lambda}(|\eta_{k}|^5 | \mathcal{F}_{k-1}) \leq c\, \epsilon^3 \Delta \langle X \rangle_{k},
\]
where $\Delta \langle X \rangle_{k}=\langle X \rangle_{k}-\langle X \rangle_{k-1}$.
From the definition of the process $V$ (see (\ref{RA-4})), it follows  that $%
\Delta \left\langle X\right\rangle _k\leq \Delta V_k=V_{k}-V_{k-1}$,
\begin{equation} \label{ee}
\mathbb{E}_{\lambda}(|\eta_{k}|^3 | \mathcal{F}_{k-1}) \leq c\, \Delta V_k \, \epsilon \quad
\textrm{and}\quad
\mathbb{E}_{\lambda}(|\eta_{k}|^5 | \mathcal{F}_{k-1}) \leq c\,  \Delta V_k\, \epsilon^3 .
\end{equation}
Thus, from (\ref{fklm59}), we obtain
\begin{eqnarray}\label{fgfgfdgk}
\mathbb{E}_{\lambda}\left( \bigg|\varphi ^{\prime \prime }\left(T_{k-1}-\frac{\vartheta_{k}\eta_{k}}{ \sqrt{A_{k}}} \right) \eta _k^3 \bigg|\mathbf{1}_{\left\{|\eta_{k}/\sqrt{A_{k}}|\leq |T_{k-1}|/2 \right\}} \Bigg| \mathcal{F}_{k-1} \right) \leq c_4\, g_{1}(T_{k-1}) \Delta V_k \, \epsilon.
\end{eqnarray}
From (\ref{fktm59}), by (\ref{ee}) and the inequality $\frac{\epsilon^2}{A_k}\geq c_{*}^{-2}$, we find
\begin{eqnarray}\label{fgfgfdgh}
\mathbb{E}_{\lambda}\left( \bigg|\varphi ^{\prime \prime }\left(T_{k-1}-\frac{\vartheta_{k}\eta_{k}}{ \sqrt{A_{k}}} \right) \eta _k^3 \bigg|\mathbf{1}_{\left\{|\eta_{k}/\sqrt{A_{k}}| > |T_{k-1}|/2 \right\}} \Bigg| \mathcal{F}_{k-1} \right) \leq   g_{2}(T_{k-1})\Delta V_k\, \epsilon,
\end{eqnarray}
where $g_{2}(t) = 2\, c (  \mathbf{1}_{\{|t|<2\}}+4 \frac{1}{t^2}  \mathbf{1}_{\{|t|\geq2\}})$.
Set $G(t)=c_4\, g_{1}(t)+g_{2}(t)$.  Then $G(t)$ is a symmetric integrable function of bounded variation, non-increasing  in $t\geq 0$. Returning to (\ref{CC-1}), by (\ref{fgfgfdgk}) and (\ref{fgfgfdgh}),  we get
\begin{eqnarray}
\left| I_1\right| \leq \mathbb{E}_{\lambda} \left[\, \sum_{k=1}^n\frac 1{6A_k^{3/2}} \mathbb{E}_{\lambda}\left(  \bigg|\varphi ^{\prime \prime }\left(T_{k-1}-\frac{\vartheta_{k}\eta_{k}}{ \sqrt{A_{k}}} \right) \eta _k^3  \bigg| \Bigg| \mathcal{F}_{k-1} \right) \right] \leq   J_1,
\end{eqnarray}
 where
\begin{equation}
J_1=c\, \epsilon \, \mathbb{E}_\lambda \sum_{k=1}^n\frac 1{A_k^{3/2}}G\left( T_{k-1} \right) \Delta V_k \,.  \label{W-0-0}
\end{equation}
Let us introduce the time change $\tau _t$ as follows:
for any real $t\in [0,T]$,
\begin{equation}
\tau _t=\min \{k\leq n:V_k>t\},\quad \textrm{where}\quad \min \emptyset =n.
\end{equation}
It is clear that, for any $t\in [0,T],$ the stopping time $\tau _t$ is
predictable. Let $(\sigma _k)_{k=1,...,n+1}$ be the increasing sequence of
moments when the increasing stepwise function $\tau _t,$ $t\in [0,T]$, has
jumps. It is clear that $\Delta V_k=\int_{[\sigma _k,\sigma _{k+1})}dt$ and
that $k=\tau _t,$ for $t\in [\sigma _k,\sigma _{k+1}).$ Since $\tau
_T=n,$ we have
\begin{eqnarray*}
\sum_{k=1}^n\frac 1{A_k^{3/2}}\, G\left( T_{k-1}\right)
\Delta V_k &=&\sum_{k=1}^n\int_{[\sigma _k,\sigma _{k+1})}\frac 1{A_{\tau
_t}^{3/2}}\, G\left( T_{ \tau_t-1}  \right) dt \\
&=&\int_0^T\frac 1{A_{\tau _t}^{3/2}}\, G\left( T_{ \tau_t-1} \right) dt.
\end{eqnarray*}
Set, for brevity, $a_t=c_{*}^2\gamma ^2+T-t.$ Since $\Delta V_{\tau _t}\leq 12\gamma ^2,$ we see that
\begin{equation}
t\leq V_{\tau _t}\leq V_{\tau _t-1}+\Delta V_{\tau _t}\leq t+12\gamma
^2,\quad t\in [0,T].  \label{BOUND-V}
\end{equation}
Taking into account that $c_{*}\geq 4,$ we have
\begin{equation}
\frac 14 \, a_t\leq A_{\tau _t}=c_{*}^2\gamma ^2+T-V_{\tau _t}\leq a_t,\quad
t\in [0,T].  \label{BOUND-A}
\end{equation}
Since $G(z)$ is symmetric and is non-increasing in $z\geq 0,$ the last bound
implies that
\begin{equation}
J_1\leq c \, \epsilon\int_0^T\frac 1{a_t^{3/2}} \, \mathbb{E}_\lambda G\left( \frac{%
u-Y_{\tau _t-1}}{a_t^{1/2}}\right) dt.  \label{W-0}
\end{equation}
By Lemma \ref{LEMMA-APX-2}, it is easy to see that
\begin{equation}
\mathbb{E}_\lambda G\left( \frac{u-Y_{\tau _t-1}}{{a_t}^{1/2}}\right) \leq
c_1\sup_z\left| \mathbb{P}_{\lambda}(Y_{\tau _t-1}\leq z)-\Phi (z)\right| +c_2\sqrt{a_t}.
\label{W-1}
\end{equation}
Since $V_{\tau _t-1}=V_{\tau _t}-\Delta V_{\tau _t}$, $ V_{\tau _t}\geq t$ (cf. (\ref{BOUND-V})) and  $\Delta V_{\tau _t}\leq
12\gamma ^2$, we get
\begin{equation}
V_n-V_{\tau _t-1}\leq V_n - V_{\tau _t}+\Delta V_{\tau _t} \leq 12\gamma ^2+T-t\leq a_t.  \label{CC-2}
\end{equation}
Thus
\begin{eqnarray*}
\mathbb{E}_\lambda \left((Y_n-Y_{\tau _t-1})^2|\mathcal{F}_{\tau _t-1} \right) &=&\mathbb{E}_{\lambda}\left(\sum_{k=\tau_t}^n\mathbb{E}_\lambda (\eta _k^2 |\mathcal{F}_{k-1})
\Bigg|\mathcal{F}_{\tau_t -1}\right) \\
&\leq &c\, \mathbb{E}_{\lambda} \left( \sum_{k=\tau_t}^n\Delta \left\langle X\right\rangle _k
\Bigg| \mathcal{F}_{\tau_t -1}\right) \\
&= &c\, \mathbb{E}_{\lambda}\left( \left\langle X\right\rangle _n-\left\langle
X\right\rangle_{\tau _t-1}|\mathcal{F}_{\tau_t -1}\right)\\
&\leq& c\, \mathbb{E}_{\lambda} \left(V_n-V_{\tau _t -1}|\mathcal{F}_{\tau_t -1}\right)    \\
&\leq& c\, a_t.
\end{eqnarray*}
Then, by Lemma \ref{LEMMA-APX-1}, we find that, for any $t\in [0,T],$%
\begin{equation}
\sup_z\left| \mathbb{P}_{\lambda}(Y_{\tau _t-1}\leq z)-\Phi (z)\right| \leq c_1 \, \sup_z\left|
\mathbb{P}_{\lambda}(Y_n\leq z)-\Phi (z)\right| +c_2\sqrt{a_t}.  \label{W-2}
\end{equation}
From (\ref{W-0}), (\ref{W-1}) and (\ref{W-2}), we obtain
\begin{equation}
J_1\leq c_1\, \epsilon \int_0^T\frac{dt}{a_t^{3/2}}\sup_z\left| \mathbb{P}_{\lambda}(Y_n\leq
z)-\Phi (z)\right| +c_2\, \epsilon \int_0^T\frac{dt}{a_t}.  \label{COMP-1}
\end{equation}
By elementary computations, we see that (since $\lambda \geq 1$)
\begin{equation}
\int_0^T\frac{dt}{a_t^{3/2}}\leq \frac c{c_{*}\lambda \epsilon }\leq
\frac c{c_{*}\epsilon }\quad \quad \mbox{and} \quad \quad \int_0^T\frac{dt}{a_t}\leq c\left| \log
  \epsilon  \right| .  \label{W-3}
\end{equation}
Then
\begin{equation}
\left| I_1\right| \leq J_1\leq \frac{c }{c_{*}}\sup_z\left| \mathbb{P}(Y_n\leq
z)-\Phi (z)\right| +c_2\, \epsilon \left| \log   \epsilon  \right|
.  \label{W-4}
\end{equation}

\emph{\textbf{b)} Control of} $I_2.$
Set $\widetilde{G}(z)=\sup_{\left| v\right| \leq 2}\psi (z+v),$ where $\psi (z)=\varphi (z)(1+z^2)^{3/2}.$
Then $\widetilde{G}(z)$ is a symmetric integrable function of bounded variation, non-increasing  in $t\geq 0$. Since $\Delta A_k=-\Delta V_k,$ we have $\left|
I_2\right| \leq I_{2,1}+I_{2,2},$ where
\begin{eqnarray*}
I_{2,1}&=&\mathbb{E}_\lambda \sum_{k=1}^n\frac 1{2A_k}\left| \varphi ^{\prime }\left(
T_{k-1}\right) \left( \Delta V_k-\Delta
\left\langle X\right\rangle _k\right) \right| ,\\
I_{2,2}&=&\mathbb{E}_\lambda \sum_{k=1}^n\frac 1{2A_k}\left| \varphi ^{\prime }\left(T_{k-1}\right) \left( \Delta \left\langle
Y\right\rangle _k-\Delta \left\langle X\right\rangle _k\right) \right| .
\end{eqnarray*}
We first deal with $I_{2,1}$. Since $\left| \varphi ^{\prime }(z)\right| \leq
\psi (z)\leq \widetilde{G}(z),$ for any real $z,$ we have
\begin{equation}
\left| \varphi ^{\prime }\left(T_{k-1}\right) \right|
\leq \widetilde{G}\left( T_{k-1}\right) .  \label{Bound-G}
\end{equation}
Note that $0\leq \Delta V_k-\Delta \left\langle X\right\rangle _k\leq
2\delta ^2\mathbf{1}_{\{k=n\}}$, $A_n = c_{*}^2\gamma ^2$ and $c_{*}\geq 4$. Then, using (\ref{Bound-G}), we get the estimations
\[
I_{2,1}\leq \frac{c_2}{c_{*}} \ \mathbb{E}_\lambda \widetilde{G}\left( T_{n-1}
\right) ,
\]
and, by (\ref{W-1}) with $G=\widetilde{G}$ and (\ref{W-2}) with $t=T,$
\[
\left| I_{2,1}\right| \leq \frac{c_1}{c_{*}}\sup_z\left| \mathbb{P}_{\lambda}(Y_n\leq z)-\Phi
(z)\right| +c_2\gamma .
\]
We next consider $I_{2,2}.$ By (\ref{f56}),  we easily obtain the bound
\[
\left| \Delta \left\langle Y\right\rangle _k-\Delta \left\langle
X\right\rangle _k\right| \leq c\lambda \epsilon \Delta \left\langle X\right\rangle _k\leq
c\lambda \epsilon \Delta V_k.
\]
With this bound, we get
\[
\left| I_{2,2}\right| \leq c\lambda \epsilon\, \mathbb{E}_\lambda \sum_{k=1}^n\frac
1{2A_k}\left|\varphi ^{\prime }\left(T_{k-1}\right)\right|
\Delta V_k.
\]
Since $\left| \varphi ^{\prime }(z)\right| \leq \psi (z) \leq \widetilde{G}(z),$ the right-hand
side can be bounded exactly in the same way as $J_1$ in (\ref{W-0-0}), with $%
A_k$ replacing $A_k^{3/2}.$ What we get is (cf. (\ref{COMP-1}))
\[
\left| I_{2,2}\right| \leq c_1\lambda \epsilon \int_0^T\frac{dt}{a_t}%
\sup_z\left| \mathbb{P}_{\lambda}(Y_n\leq z)-\Phi (z)\right| +c_2\lambda \epsilon \int_0^T%
\frac{dt}{a_t^{1/2}}.
\]
By elementary computations, we see that
\[
\int_0^T\frac{dt}{a_t^{1/2}}\leq \int_0^T\frac{dt}{\sqrt{T-t}}\leq c_2,
\]
and, taking into account that $a_t\geq c_{*}^2\gamma ^2,$
\[
\int_0^T\frac{dt}{a_t}\leq \frac{c_1}{c_{*}\lambda \epsilon }\int_0^T%
\frac{dt}{a_t^{1/2}}\leq \frac{c_2}{c_{*}\lambda \epsilon }.
\]
Then
\[
\left| I_{2,2}\right| \leq \frac{c_1}{c_{*}}\sup_z\left| \mathbb{P}_{\lambda}(Y_n\leq
z)-\Phi (z)\right| +c_2\lambda \epsilon .
\]
Collecting the bounds for $I_{2,1}$ and $I_{2,2},$ we get
\begin{equation}
\left| I_2\right| \leq \frac{c_1}{c_{*}}\sup_z\left| \mathbb{P}_{\lambda}(Y_n\leq z)-\Phi
(z)\right| +c_2\gamma .  \label{F-2}
\end{equation}

\emph{\textbf{c)} Control of} $I_3.$  By Taylor's expansion,
\[
I_3=\frac{1}{8}\,\mathbb{E}_\lambda \sum_{k=1}^n\frac 1{(A_k-\vartheta _k\Delta A_k)^2}\varphi
^{\prime \prime \prime }\left( \frac{u-Y_{k-1}}{\sqrt{A_k-\vartheta _k\Delta A_k} }
\right) \Delta A_k^2.
\]
Since $\left| \Delta A_k\right| =  \Delta V_k  \leq 12\gamma ^2$
and $c_{*}\geq 4,$ we have
\begin{equation}
A_k\leq A_k-\vartheta _k\Delta A_k\leq c_{*}^2\gamma ^2+T-V_k+12\gamma ^2\leq
2A_k. \label{eqsfion}
\end{equation}
Using (\ref{eqsfion}) and the inequalities $\left| \varphi ^{\prime \prime \prime
}(z)\right| \leq \psi (z) \leq \widetilde{G}(z)$, we obtain
\[
\left| I_3\right| \leq c\gamma ^2\mathbb{E}_\lambda \sum_{k=1}^n\frac 1{A_k^2} \ \widetilde{G}\left(
\frac{T_{k-1}}{\sqrt{2 }}\right) \Delta V_k.
\]
Proceeding in the same way as for estimating $J_1$ in (\ref{W-0-0}), we get
\begin{equation}
\left| I_3\right| \leq \frac{c_1}{c_{*}}\sup_z\left| \mathbb{P}_{\lambda}(Y_n\leq z)-\Phi
(z)\right| +c_2\gamma .  \label{F-3}
\end{equation}

We are now in a position to end the proof of Lemma \ref{LEMMA4}. From  (\ref{lla}),  using (\ref{W-4}), (\ref
{F-2}) and (\ref{F-3}), we find
\begin{eqnarray*}
 \left| \mathbb{E}_\lambda \Phi _u(Y_n,A_n)-\mathbb{E}_\lambda \Phi _u(Y_0,A_0)\right|   \leq \frac{c_1}{c_{*}}\sup_z\left| \mathbb{P}_{\lambda}(Y_n\leq z)-\Phi
(z)\right| +c_2(\lambda \epsilon +\epsilon \left| \log
\epsilon  \right| +\delta ).
\end{eqnarray*}
Implementing the last bound in (\ref{llaa}), we come to
\[
\sup_z\left| \mathbb{P}_{\lambda}(Y_n\leq z)-\Phi (z)\right| \leq \frac{c_1}{c_{*}}\sup_z\left|
\mathbb{P}_{\lambda}(Y_n\leq z)-\Phi (z)\right| +c_2(\lambda \epsilon +\epsilon \left|
\log   \epsilon  \right| +\delta ),
\]
from which, choosing $c_{*}=\max\{ 2c_1, 4\}$, we get
\begin{equation}
\sup_z\left| \mathbb{P}_{\lambda}(Y_n\leq z)-\Phi (z)\right| \leq 2c_2(\lambda \epsilon
+\epsilon \left| \log  \epsilon  \right| +\delta ),
\label{F-4}
\end{equation}
which proves Lemma \ref{LEMMA4} for $1\leq \lambda < \epsilon ^{-1}$.

For $0\leq \lambda <1$, we can prove Lemma \ref{LEMMA4} similarly by taking $\gamma = \epsilon|\log \epsilon| +\delta $.  We only need to note that in this case, instead of (\ref{W-3}),
\begin{equation}
\int_0^T\frac{dt}{a_t^{3/2}} \leq
\frac c{c_{*}\epsilon|\ln \epsilon| } \quad \quad \mbox{and} \quad \quad  \int_0^T\frac{dt}{a_t}\leq c\left| \log
 \epsilon \right|.
\end{equation}

\section{Equivalent conditions} \label{seca}
In the following we give several equivalent conditions to the Bernstein condition (A1$'$). In the independent case equivalent conditions
can be found in Saulis and Statulevi\v{c}ius \cite{SS78}. For the convenience of the readers and motivated by the fact that in \cite{SS78} the conditions
are rather different from those used here, we decided to include independent proofs.
\begin{proposition}\label{s}
The following three conditions are equivalent:\\
\emph{(I)} Bernstein's condition (A1$'$). \\
\emph{(II) (Sakhanenko's condition)} There exists some positive absolute  constant $K$ such that
  \[
  K\,\mathbb{E}(|\eta_{i}|^3\exp\{K|\eta_{i}|\}|\mathcal{F}_{i-1})\leq \mathbb{E}(\eta_{i}^2|\mathcal{F}_{i-1}),\ \ \ \  \textrm{for} \  1\leq i\leq n.
  \]
\emph{(III)} There exists some positive absolute  constant $ \rho$ such that
\[
\mathbb{E}(|\eta_{i}| ^{k}  | \mathcal{F}_{i-1})\leq \frac12 k!\rho^{k-2} \mathbb{E}(\eta_{i}^2 | \mathcal{F}_{i-1}), \ \ \ \textrm{for}\ k\geq 3\ \ \textrm{and}\ \ 1\leq i\leq n.
\]
\end{proposition}

\noindent\textbf{Proof.} First we prove that (I) implies (II). Let $t \in (0, 1)$.
By condition (I) and Lemma \ref{l11}, we find that
\begin{eqnarray}
 E(|\eta_{i}|^3e^{t H^{-1}|\eta_{i}|} | \mathcal{F}_{i-1}) &=& \sum_{k=0}^{\infty} \frac{ (t H^{-1})^{k}}{k!}  E(|\eta_{i}|^{k+3} | \mathcal{F}_{i-1}) \nonumber \\
  &\leq& \sum_{k=0}^{\infty} \frac{ (t H^{-1})^{k}}{k!}(k+3)! H^{k+1}E(\eta^2_{i} | \mathcal{F}_{i-1}) \nonumber \\
  &\leq& H E(\eta^2_{i} | \mathcal{F}_{i-1})\sum_{k=0}^{\infty} \frac{ (k+3)!}{ k!} \, t^k \nonumber \\
  &=:&f(t) H E(\eta^2_{i} | \mathcal{F}_{i-1}). \label{hndsf}
\end{eqnarray}
Since $g(t)= t f(t)$ is a continuous function in $[0, \frac12]$ and satisfies $g(0)=0$ and $  g(\frac12) \geq 3$, there exists $t_0 \in (0, \frac12)$ such that $g(t_0)=1$. Taking $K=t_0 H^{-1}$, we obtain condition (II) from (\ref{hndsf}).

Next we show that (II) implies (III).  By the elementary inequality $x^k \leq k! \, e^x$, for $k \geq0$ and $x \geq 0$,
it follows that, for $k\geq 3$,
\begin{eqnarray*}
 \mathbb{E}(|\eta_{i}| ^{k}  | \mathcal{F}_{i-1}) &= &  \mathbb{E}( |\eta_{i}| ^{3}K^{3-k} |K \eta_{i}| ^{k-3}  | \mathcal{F}_{i-1}) \\
 &\leq& (k-3)! \, K^{3-k} \mathbb{E}( |\eta_{i}| ^{3} \exp\{|K \eta_{i}|\}   | \mathcal{F}_{i-1}) .
\end{eqnarray*}
Using condition (II),  for $k\geq 3$,
\begin{eqnarray*}
 \mathbb{E}(|\eta_{i}| ^{k}  | \mathcal{F}_{i-1})
 \leq  (k-3)! \, K^{2-k} \mathbb{E}( \eta_{i} ^2  | \mathcal{F}_{i-1}) \leq \frac12 k! \rho^{k-2}  \mathbb{E}(\eta_{i}^2 | \mathcal{F}_{i-1}),
\end{eqnarray*}
where $\rho =  \frac1K $, which proves condition (III).

It is obvious that (III) implies (I) with $H=\rho$. \hfill\qed

\begin{proposition} If $\eta_1,...,\eta_n$ are i.i.d., then Bernstein's condition, Cram\'{e}r's condition and  Sakhanenko's  condition are all equivalent.
\end{proposition}

\noindent\textbf{Proof.} According to Theorem \ref{s}, we only need to prove that Cram\'{e}r's condition and Bernstein's condition are equivalent.
We can assume that, a.s., $\eta_{1} \ne  0.$

First, from (\ref{sbopdf}), we find that Bernstein's condition (A1$'$) implies Cram\'{e}r's condition:
\begin{eqnarray*}
\mathbb{E} e^{ \frac12 H^{-1} \eta_{1}}
  < \infty.
\end{eqnarray*}

Second, we show that Cram\'{e}r's condition, i.e. $\mathbb{E} e^{c_0^{-1}|\eta_{1}|}:=c_1< \infty$, implies Bernstein's condition (A1$'$).
 By the inequality $x^k \leq k! \, e^x$, for $k \geq0$ and $x \geq 0$, it  follows that
\begin{eqnarray*}
|\mathbb{E} \eta_{1}^k|  \leq c_0^k \, \mathbb{E} |c_0^{-1} \eta_{1}|^k \leq k! \, c_0^{k} \, \mathbb{E} e^{c_0^{-1} |\eta_{1}|}
& =& k! \, c_0^{k} \, c_1.
\end{eqnarray*}
Then, it is easy to see that, for $k\geq 3$,
\begin{eqnarray*}
|\mathbb{E} \eta_{1}^k|  \leq  \frac12 \,  k! \, c_0^{k-2} \frac{2c_0^{2}c_1}{\sigma^2} \,  \mathbb{E} \eta_{1}^2 \leq  \frac12 \,  k! \, H^{k-2}  \,  \mathbb{E} \eta_{1}^2,
\end{eqnarray*}
where $\sigma^2=\mathbb{E} \eta_{1}^2$ and $H=\max\bigg\{c_0,  \frac{2c_0^{3}c_1}{\sigma^2} \bigg\},$ which proves that condition (A1$'$) is satisfied.
\hfill\qed

\section*{Acknowledgements}
We would like to thank the two referees for their helpful remarks and suggestions. The work
has been partially supported by the National Natural Science Foundation of China (Grant no.\ 11101039 and
 Grant no.\ 11171044), and Hunan Provincial Natural Science Foundation of China (Grant No.11JJ2001).
Fan was partially supported by the Post-graduate Study Abroad Program sponsored by China Scholarship Council.

\end{document}